\documentclass[12pt]{amsart}
\usepackage{tikz, stmaryrd}
\usetikzlibrary{matrix,arrows,decorations.pathmorphing}
\usepackage[norelsize]{algorithm2e}
\usepackage{color}
\usepackage{amsmath,amssymb, amsfonts, amscd}
\usepackage[all]{xy}
\usepackage{verbatim}
\usepackage{enumerate}
\usepackage{yhmath}
\usepackage{xcolor}

\newtheorem{thm}{Theorem}[section]
\newtheorem{cor}[thm]{Corollary}
\newtheorem{lem}[thm]{Lemma}
\newtheorem{obs}[thm]{Observation}

\newtheorem{prop}[thm]{Proposition}

\theoremstyle{definition}

\newtheorem{defn}[thm]{Definition}
\newtheorem{notation}[thm]{Notation}

\newtheorem{rem}[thm]{Remark}
\numberwithin{equation}{section}

\DeclareFontFamily{U}{rsf}{} \DeclareFontShape{U}{rsf}{m}{n}{
  <5> <6> rsfs5 <7> <8> <9> rsfs7 <10->  rsfs10}{}
\DeclareMathAlphabet{\mathscr}{U}{rsf}{m}{n}

\newcommand{\Oc}{\mathcal O}
\newcommand{\A}{A}
\newcommand{\B}{B}

\newcommand{\D}{D}

\newcommand{\Zb}{{\mathbb Z}}

\renewcommand{\imath}{\sqrt{-1}}

\newcommand{\Hc}{\mathcal H}

\newcommand{\Ab}{\mathbb A}

\newcommand{\Nb}{\mathbb N}

\newcommand{\Lb}{\mathbb L}
\newcommand{\Rb}{\mathbb R}
\newcommand{\Cb}{\mathbb C}
\newcommand{\Qb}{\mathbb Q}

\DeclareMathOperator{\Hom}{Hom}
\DeclareMathOperator{\colim}{colim}
\DeclareMathOperator{\spec}{spec}
\DeclareMathOperator{\id}{e}
\DeclareMathOperator{\coim}{coim}
\DeclareMathOperator{\im}{im}
\DeclareMathOperator{\coker}{coker}

\newcommand{\uHom}{\underline{\Hom}}
\newcommand{\ootimes}{\overline{\otimes}}
\newcommand{\wotimes}{\widehat{\otimes}}

\newcommand{\tC}{\text{\bfseries\sf{C}}}
\newcommand{\ttsC}{\text{\bfseries\sf{sC}}}

\newcommand{\mM}{\mathcal{M}}

\newcommand{\mW}{\mathcal{W}}

\newcommand{\mT}{\mathcal{T}}
\newcommand{\mO}{\mathcal{O}}
\newcommand{\mB}{\mathcal{B}}

\newcommand{\ttMod}{{\text{\bfseries\sf{Mod}}}}
\newcommand{\ttBorn}{{\text{\bfseries\sf{Born}}}}
\newcommand{\ttSNrm}{{\text{\bfseries\sf{SNrm}}}}
\newcommand{\ttNrm}{{\text{\bfseries\sf{Nrm}}}}
\newcommand{\ttCBorn}{{\text{\bfseries\sf{CBorn}}}}
\newcommand{\ttSBorn}{{\text{\bfseries\sf{SBorn}}}}

\newcommand{\ttBan}{{\text{\bfseries\sf{Ban}}}}

\newcommand{\ttVect}{{\text{\bfseries\sf{Vect}}}}
\newcommand{\ttComm}{{\text{\bfseries\sf{Comm}}}}

\newcommand{\ttC}{{\tC}}

\newcommand{\ttAff}{{\text{\bfseries\sf{Aff}}}}
\newcommand{\ttAfnd}{{\text{\bfseries\sf{Afnd}}}}

\newcommand{\ttSets}{{\text{\bfseries\sf{Sets}}}}

\newcommand{\ttInd}{{\text{\bfseries\sf{Ind}}}}

\DeclareMathOperator{\sep}{sep}
\newcommand{\ie}{\emph{i.e.}    }
\newcommand{\la}{\lambda}
\renewcommand{\a}{\alpha}
\newcommand{\be}{\beta}
\def\then{\Rightarrow}
\newcommand{\ol}{\overline}

\newcommand{\what}{\widehat}
\def\rhook{{\hookrightarrow}}
\def\lt{\langle}
\def\gt{\rangle} 
\def\diss{{\rm diss\,}}

\def\sup{{\rm sup\,}}

\def\limpro{\mathop{\lim\limits_{\displaystyle\leftarrow}}}
\def\limind{\mathop{\lim\limits_{\displaystyle\rightarrow}}}

\numberwithin{equation}{section}

\begin{document}

\title[]{Dagger Geometry As Banach Algebraic Geometry}
\author{Federico Bambozzi, Oren Ben-Bassat}\thanks{}%
\address{Federico Bambozzi, Fakult\"{a}t f\"{u}r Mathematik,
Universit\"{a}t Regensburg,
93040 Regensburg, 
Germany}
\bigskip
\email{Federico.Bambozzi@mathematik.uni-regensburg.de}
\address{Oren Ben-Bassat, Department of Mathematics, University of Haifa, Haifa, Israel}
\hfill \break
\bigskip
\email{ben-bassat@math.haifa.ac.il}
\dedicatory{}
\thanks{
        The first author acknowledges the support of the University of Padova by MIUR PRIN2010-11 "Arithmetic Algebraic Geometry and Number Theory", the University of Regensburg and with the support of the DFG funded CRC 1085 "Higher Invariants. Interactions between Arithmetic Geometry and Global Analysis" that permitted him to work on this project. The second author acknowledges the University of Oxford and the support of the European Commission under the Marie Curie Programme for the IEF grant which enabled this research to take place. The contents of this article reflect the views of the two authors and not the views of the European Commission. We would like to thank Francesco Baldassarri, Elmar Grosse-Kl\"onne, Kobi Kremnizer, Fr\'{e}d\'{e}ric Paugam, and J\'{e}r\^{o}me Poineau for interesting conversations.}

\begin{abstract}
\noindent 
In this article, we apply the approach of relative algebraic geometry towards analytic geometry to the category of bornological and Ind-Banach spaces (non-Archimedean or not). We are able to recast the theory of Grosse-Kl\"onne dagger affinoid domains with their weak G-topology in this new language. We prove an abstract recognition principle for the generators of their standard topology (the morphisms appearing in the covers). We end with a sketch of an emerging theory of dagger affinoid spaces over the integers, or any Banach ring, where we can see the Archimedean and non-Archimedean worlds coming together.
\end{abstract}
\maketitle 
\tableofcontents
\section{Introduction} 

Analytic geometry tends to lack many of the tools and the clear presentation of algebraic geometry. One example is the theory of quasi-coherent sheaves in analytic geometry which has experienced a plethora of different definitions (see for instance the discussion in \cite{BBP}). There is also not an obvious way of dealing with analytic spaces which are not finite dimensional, or in general with various limits and colimits. One sees many objects described without satisfying universal properties. In short, there is no Grothendieck style approach to analytic geometry yet established. On the other hand, many infinite dimensional constructions in algebriac geometry could be interpreted as constructions in analytic geometry that are in some sense, finite dimensional. Due to historical reasons, p-adic and complex analytic geometry also developed separately and lack a common language.  We feel like the correct language to use is that of relative algebraic geometry. This means that we do geometry relative to a closed (unital) symmetric monoidal category in the sense of \cite{TV}, \cite{TVe}, \cite{TVe2}. To\"{e}n and Vezzosi  define formal Zariski open immersions between affine schemes relative to the closed symmetric monoidal quasi-abelian category satisfying some mild extra conditions. The focus of the current article is understanding the meaning of these conditions in cases relevant to our new approach to analytic geometry.  This idea of this approach is due to Kobi Kremnizer and this article is one in a series of papers including \cite{BeKr}, \cite{BeKr2} and others to appear which will develop it. The main symmetric monoidal categories that are used in this article are the category of complete Bornological vector spaces of convex type over a complete, non-trivially valued field and a generalization of this: the categories of Ind-Banach modules over a complete normed ring. Some of the most interesting examples are the non-Archimedean rings $\mathbb{Q}_{p}$ and $\mathbb{C}((t))$, the Novikov ring, and the Archimedean rings, $\mathbb{Z}, \mathbb{R}, \mathbb{C}$. Some of these already have their own analytic geometry associated with it and we give a unified description.
In the present article, we begin by showing that the theory of dagger affinoid subdomains of dagger affinoid spaces has an abstract description in terms of derived categories (or derived tensor products).  The morphisms appearing in the covers defining the topology on this category are characterized as being homotopy monomorophisms. The first author showed in \cite{Bam} that bornological algebras are a natural setting for dagger analytic geometry, or the geometry of over-convergent analytic functions. Recently, similar algebras were also studied in relationship to Wiener theorems by Alpay and Salomon \cite{AlSo}.  It was shown by Grosse-Kl\"onne \cite{GK} that algebras of over-convergent analytic functions are needed to give a reasonable geometric realization of de Rham cohomology in non-Archimedean geometry.  

The main technical tool we use is the correct version of the derived category, which comes from the quasi-abelian setting \cite{SchneidersQA}.  In order to show that we can recover well known notions in complex and $p$-adic analytic geometry, we first look at a special case of our framework: we limit the Banach ring $R$ to a non-trivially valued field and we consider the category of algebras of a certain special form. The categories we consider in this article is the opposite of the category of dagger affinoid algebras over a valuation field. We examine the abstractly defined topology (the formal Zariski topology coming from the derived setting of To\"{e}n-Vezzosi) and show that it agrees with the standard topology. Our main results are Theorems \ref{thm:LocToHoEp} and \ref{thm:HoEpToEmb}.

An advantage of this foundational approach will be towards a theory of analytic stacks related to recent work (\cite{Br}, \cite{PY}, \cite{Ul}) , analytic D-modules (as in \cite{AW}), constructible sheaves, Arakelov Theory and more. In future work, we will look at geometry over the Adeles in this context and develop a form of descent which can prove results for schemes over the analytic integers by proving them over $\mathbb{R}$ and the $\mathbb{Q}_{p}$. We also think to be able to give a unified treatment of many theorems from algebraic and analytic geometry including Grauert's pushforward theorems, theorem A and B, and others. 

In an upcoming work, we will extend the results in this article to the Stein setting (increasing unions of affinoids or dagger affinoids) as well as the quasi-Stein setting. This allows one to use relative algebraic geometry to understand spaces without boundary, as well as ones with partial boundary. In both cases we can characterize the morphisms between them which are usually known as open immersions.  

\begin{notation}
In this article the following notation will be used, if not otherwise stated:
\begin{itemize}
 \item $k$ will denote a field complete with respect to a fixed absolute value, Archimedean or non-Archimedean.
 \item we will denote with $k^\circ = \{ x \in k | |x| \le 1 \}$, the subset of power-bounded elements of $k$.
 \item $\ttVect_{k}$ is the closed symmetric monoidal category of vector spaces (with no extra structure) over a field $k$.
  \item $R$ will denote a complete normed ring, the norm can be Archimedean or non-Archimedean.
 \item If $\ttC$ is a category we will use the notation $X \in \ttC$ to mean $X$ is an object of $\ttC$.
 \item $\ttSNrm_R$ the category of semi-normed modules over $R$, remarking that, if not otherwise stated, by a semi-normed space over non-Archimedean base field we mean a $R$-vector space equipped with a non-Archimedean semi-norm.
 \item $\ttNrm_R$ the category of normed modules over $R$.
 \item $\ttBan_R$ the category of Banach modules over $R$.
 \item For $V\in \ttSNrm_R$, $V^{s}= V/\overline{(0)}\in \ttNrm_{R}$ is the seperation and $\widehat{V} \in \ttBan_{k}$ is the seperated completion.
 \item $\ttBorn_{k}$ the category of bornological vector spaces of convex type over $k$.
 \item For $E \in \ttBorn_{k}$ and $B$ a bounded absolutely convex subset, (a bounded disk) of $E$  then $E_B$ is the vector subspace of $E$ spanned by elements of $B$  and equipped with the gauge semi-norm (also called the Minkowski functional) defined by $B$.
 \item  For $E \in \ttBorn_{k}$, $\mathcal{B}_{E}$ denotes the category of bounded, absolutely convex subsets of $E$.
  \item  For $E \in \ttBorn_{k}$, $\mathcal{B}^{c}_{E}$ denotes the category of bounded, absolutely convex subsets $B$ of $E$ for which $E_{B}\in \ttBan_{k}$ 
 \item  If $\ttC$ is a category then $\widehat{\ttC}$ denote the category of contravariant functors from $\ttC$ to the category of sets with morphisms being natural transformations.
  \item If $\ttC$ is a category then $\ttInd(\ttC)$ will denote the category of Ind-systems over $\ttC$.
 \item If $V_i$ are objects in a concrete category, $\underset{i\in I}\times V_{i}$ denotes the product of the underlying sets of the $V_i$.
 \item If $M$ is some $R$-module with extra structure for some ring $R$ then $M^{\times}$ refers to the set $M-\{0\}$;
 \item $\bigoplus_{i\in I} M_i$ is the direct sum a collection of modules $\{M_i\}_{i\in I}$ for some ring which is the coproduct in the category of (algebraic) modules.
 \item The notation $\limind$ refers to a colimit (also known as inductive or direct limit) of some direct system in a category.
  \item the notation $\limpro$ refers to a limit (also known as projective or inverse limit) of some direct system in a category.
  \item We will often consider polyradii $\rho = (\rho_i) \in \Rb_+^n$ and follow the convention that the notation $\rho < \rho'$ means that $\rho$ and $\rho'$ have the same number of components and every component of $\rho$ is strictly less than the correspondent component of $\rho'$. 
\end{itemize}
\end{notation}

\section{Background Material}

\subsection{Ind Categories}

Here we collect some basic results in the theory of Ind-categories and fix the notation that will be used.

\begin{defn}
Let $\ttC$ be a category and $I$ a small filtered category. An \emph{Ind-object} of $\ttC$ is a functor $X: I \to \ttC$, denoted by $X = (X_i)_{i \in I}$. 
\end{defn}

Let's denote by $h: \ttC \to \widehat{\ttC} $ the Yoneda embedding. Then, for any Ind-object one obtains a pre-sheaf
\[ L_I(X)  = \limind_{i \in I} h_{X_i}. \]
This object will be often denoted with  ${``\limind"}_{i\in I} X_i.$
An object $F \in \widehat{\ttC}$ which is isomorphic to $L_I(X)$, for some Ind-object $X$ is called \emph{Ind-representable pre-sheaf}. 

\begin{defn}
Let $X = (X_i)_{i \in I}, Y = (Y_j)_{j \in J}$ be two Ind-objects of $\ttC$, then we define
\[ \Hom(X, Y) = \Hom_{\widehat{\ttC}}(L_I(X), L_J(Y)). \]
With this definition of morphisms the class of Ind-objects of $\ttC$ forms a category which is called the \emph{category of Ind-objects} of $\ttC$ (or simply the Ind-category) and denoted by $\ttInd(\ttC)$. 
\end{defn}

Explicitly, one finds that given $X = (X_i)_{i \in I}$ and $Y = (Y_j)_{j \in J}$
\begin{equation}
\label{eqn:IndHom} \Hom(X, Y) = \limpro_{i \in I} \limind_{j \in J}\Hom(X_i, Y_j). 
\end{equation}
Hence the Yoneda embedding $\ttC \to \widehat{\ttC}$ factors through $\ttInd(\ttC)$, and we denote \begin{equation} \label{eqn:embedding}
 \iota: \ttC \to \ttInd(\ttC)
 \end{equation}
 the fully faithful embedding of categories.

\begin{rem}
With this definition of morphisms one can see, writing down explicitly the commutative diagrams involved, that for any given morphism of Ind-objects $f: {``\limind"}_{i \in I} X_i \to {``\limind"}_{j \in J} Y_j$ there exists a small filtered category $K$ and functors $X':K\to\ttC$ and $Y':K \to \ttC$ and a natural transformation of functors $N: X' \to Y'$ such that $L_I(X) \cong L_K(X')$, $L_J(Y) = L_K(Y')$ and $L(N) =f$. 
The details of this operation are explained in Meyer's book \cite{M}, page 54 remark 1.133, and in appendix 3 of \cite{ArMaz}.  We will frequently use this operation of \emph{re-indexing} a morphism.
\end{rem}

The operation of re-indexing can not be defined in a functorial way for any diagram of maps.

%

\begin{defn} \label{defn:level_representation}
Let $D: I \to \ttInd(\ttC)$ be a diagram in $\ttInd(\ttC)$. A \emph{level representation} for $D$ is the data of a filtered category $K$ and a functor $X: I \times K \to \ttC$ such that
\begin{enumerate}
\item for any $i \in I$ there exist an isomorphism ${``\limind"}_{k \in K} X_{i,k} \cong D_i$;
\item for every map $\phi: i \to j$ in $I$ the morphism $X_\phi: {``\limind"}_{k \in K} X_{i,k} \to {``\limind"}_{k \in K} X_{j,k}$ and the isomorphisms in (1) fit into a commutative diagram
\[\xymatrix{``\limind"_{k \in K} X_{i,k} \ar[d] \ar[r] & ``\limind"_{k \in K} X_{j,k} \ar[d]  \\ D_{i}  \ar[r]_{D(\phi)} & D_{j}   }.
\]

\end{enumerate}
\end{defn}

\begin{prop} \label{prop:reindexing}
Let $\Delta$ be a loopless finite category and let $D: \Delta \to \ttInd(\ttC)$ be a diagram of type $\Delta$. Then, $D$ admits a level representation.
\end{prop}
{\bf Proof.}

This proposition is proved for the case of the category of Pro-objects of $\ttC$ in \cite{ArMaz}, Proposition 3.3 of appendix. In it not difficult to see that the same argument works for Ind-objects.

\hfill $\Box$


%
%


We finish this section by recalling some well known results on calculation of limits and colimits in the category $\ttInd(\ttC)$. To do so, we introduce the functor $L: \ttInd(\ttC) \to \what{\ttC}$ which is defined 
\[ L(X) = L_I(X) \]
for any Ind-object $X = (X_i)_{i \in I}$.

\begin{lem}
The functor $L: \ttInd(\ttC) \to \what{\ttC}$ commutes with all (projective) limits and filtered colimits.
\end{lem}

{\bf Proof.}
The first assertion is in Proposition 8.9.1 of the first expos\'e of \cite{SGA4}, while the second one is in Proposition 8.5.1 of the same expos\'e.

\hfill $\Box$

\begin{lem} 
Given a small filtered category $I$ and a functor $I \to \ttC$, then the corresponding object of $\ttInd(\ttC)$ is isomorphic to the colimit in $\ttInd(\ttC)$ of the diagram obtained by the composition $I \to \ttC \to \ttInd(\ttC)$.
\end{lem}

{\bf Proof.}

This follows immediately from checking that ${``\limind"}_{i\in I} X_i$ satisfies the universal property of the colimit  ${\limind}_{i\in I} \iota(X_i)$ taken in $\ttInd(\ttC).$ That is immediate from Equation (\ref{eqn:IndHom}).

\hfill $\Box$

\begin{lem} \label{lemma:ind_limits}
If $\ttC$ has all finite limits and finite colimits then the category $\ttInd(\ttC)$ has all limits and colimits.
\end{lem}
{\bf Proof.}

By \cite{SGA4} Proposition 8.9.1, expos\'e 1, $\ttInd(\ttC)$ has all limits. A construction of products and coproducts in $\ttInd(\ttC)$ for any category $\ttC$ with finite products and finite coproducts can be found in \cite{M}, page 55. General colimits in $\ttInd(\ttC)$ can be formed from coproducts and cokernels and cokernels can be constructed objectwise.
\hfill $\Box$

An important special case is:

\begin{lem} \label{lem:coproducts}
Suppose that $\ttC$ is a category admitting finite coproducts. 
The coproduct in $\ttInd(\ttC)$ of a collection of objects $X_{i} \in \ttC$ indexed by $i \in I$ is given by 
\[{``\limind"}_{J \in Fin(I)} \coprod_{j \in J} X_{j}.
\]
where $Fin(I)$ is the small filtered category of finite subsets $J \subset I$.
\end{lem}

\subsection{Quasi-abelian categories}

In this section we recall some definitions and results from the theory of quasi-abelian categories, mainly from \cite{SchneidersQA}.

\begin{defn}\label{defn:StrictExact} A sequence of morphisms \[ E' \stackrel{u}\longrightarrow E \stackrel{v}\longrightarrow E''
        \]
        is called 
        \begin{itemize}
                \item \emph{strictly exact at $E$} if $u$ is strict, $v\circ u=0$ and the natural morphism $\im(u) \to \ker(v)$ is an isomorphism;
                \item \emph{strictly coexact at $E$} if $v$ is strict, $v\circ u=0$ and the natural morphism $\im(u) \to \ker(v)$ is an isomorphism.
        \end{itemize}
        A long sequence of morphisms is called \emph{strictly exact} (resp. \emph{strictly coexact}) if it is \emph{strictly exact} (resp. \emph{strictly coexact}) at any term.
\end{defn}

We remark that for a short exact sequence the concepts of strict exactness and strict coexactness coincides.

\begin{defn}An additive functor $F$ between quasi-abelian categories is called \emph{exact} if it transforms any strictly (co)exact sequence
        \[0 \to E'\to E\to E''\to 0
        \]
        into the strictly (co)exact sequence 
        \[0 \to F(E')\to F(E)\to F(E'')\to 0.
        \]
\end{defn}
For the sake of clarity we include a precise description of strict morphisms, monomorphisms, strict monomorphisms, epimorphisms and strict epimorphisms in the Ind-category of a quasi-abelian category.

\begin{prop} \label{prop:strict_morphism_ind}
        Let $\ttC$ be a quasi-abelian category and $f: M  \to N $ be a morphism in $\ttInd(\ttC)$, then
        \begin{enumerate}
                \item $f$ is a monomorphism (resp. epimorphism) if and only if there exists a re-indexing $\limind f_k: {``\limind"}_{k \in K} M_k \to ``\limind" N_k$  of $f$ such that $f_k$ is a monomorphism (resp. epimorphism); 
                \item $f$ is a strict morphism if and only if there exists a re-indexing   $\limind f_k: {``\limind"}_{k \in K} M_k \to {``\limind"}_{k \in K} N_k$  of $f$ such that $f_k$ is a strict morphism;               
                \item $f$ is a strict monomorphism (resp. strict epimorphism) if and only if there exists a re-indexing  $\limind f_k: {``\limind"}_{k \in K} M_k \to {``\limind"}_{k \in K} N_k$  of $f$ such that $f_k$ is a strict monomorphism (resp. strict epimorphism).
        \end{enumerate}
\end{prop}
{\bf Proof.}

The last claim follows from the combination of the first two. $f$ is a monomorphism if and only if $\ker f \cong 0$, which means that $\ker f$ is isomorphic to the constant object $0$ in $\ttInd(\ttC)$. Hence there exists a re-indexing such that $\ker f \cong (\ker f_k)_{k \in K} = (0)_{k \in K}$ and for this re-indexing $f_k$ is a monomorphism for any $k$. The same reasoning works for epimorphisms and cokernels.

To prove that claim on strict morphisms is enough to notice that a strict morphism is by definition a morphism such that the canonical morphism $\coim(f) \to \im(f)$ is an isomorphism and $\im f = \ker(N \to \coker f)$ and $\coim f = \coker(\ker f \to M)$. So, since the kernels and cokernels are calculated term by term we obtain that there exists a re-indexing of $M$ and $N$, $\limind f_k: ``\limind" M_k \to ``\limind" N_k$, such that $f_k$ is a strict morphism for any $k$.

\ \hfill $\Box$
\begin{defn}Let $\ttC$ be an additive category with finite limits and colimits. We call an object $P$ \emph{projective} if for all strict epimorphisms $V\to W$, the corresponding map $\Hom(P, V) \to \Hom(P,W)$ is surjective.
\end{defn}
\begin{lem}\label{lem:FiltProj} Let $\ttC$ be an additive category with finite limits and colimits. Any filtered colimit of projectives is projective.
\end{lem}
{\bf Proof.} Let $V\to W$ be a strict epimorphism. Let $P = \limind_{i\in I} P_i$ be a  filtered colimit of projectives. Then $\Hom(P, V)\to \Hom(P,W)$ is a cofiltered limit of surjective maps of sets $\Hom(P_i, V)\to \Hom(P_i,W)$. Hence it is surjective by the classical Mittag-Leffler lemma.
\ \hfill $\Box$
\begin{defn}
	A quasi-abelian category $\ttC$ is said to have \emph{enough projectives} if for any object $X \in \ttC$ there exists a projective object $P \in \ttC$ and a strict epimorphism $P \to X$.
\end{defn}

If the category $\ttC$ is quasi-abelian and has enough projective objects then the same holds for $\ttInd(\ttC)$, which is also a quasi-abelian category. Indeed, if  $\mathcal{E}= ``\colim"_{i \in I} E_{i}$ is an element of $\ttInd(\ttC)$ and if we have strict epimorphisms $\kappa_{E_i}: P(E_i) \to E_i$ where $P(E_i)$ are projective then define 
\[P(\mathcal{E}) = \coprod_{i \in I} P(E_{i})
\]
where the coproduct is taken in $\ttInd(\ttC).$  Then $P(\mathcal{E})$ is projective in $\ttInd(\ttC)$ and using the strict epimorphisms $\kappa_{E_i}$ one gets a strict epimorphism $\kappa_{\mathcal{E}}:P(\mathcal{E})\to \mathcal{E}$.

\begin{defn}\label{def:generator} (Definition 5.2.1 of \cite{KS})
        A generator in a category $\ttC$ is an object $G$ of $\ttC$ such that the functor 
        \[\ttC^{op} \to \ttSets
        \]
        given on objects by
        \[V \mapsto \Hom(V, G)
        \]
        has the property that any morphism $f: V\to W$ which induces an isomorphism 
        \[\Hom(G,f): \Hom(G,V)\to \Hom(G,W)\] is an isomorphism.
\end{defn}

\begin{defn}\cite{SchneidersQA} A quasi-abelian category is called elementary if it is cocomplete and has a small, strictly generating set of tiny projective objects.
\end{defn}

\begin{lem}\label{lem:hasgen}Any elementary quasi-abelian category has a generator in the sense of Definition \ref{def:generator}.
\end{lem}
{\bf Proof.}
Let $\ttC$ be elementary quasi-abelian category. Let $\{P_{s}|s\in S\}$ be a small, strictly generating set of tiny projective objects of $\ttC$. Let $G= \coprod_{s\in S} P_{s}$. Let  $f: V\to W$ be a morphism which induces an isomorphism $\Hom(G,f): \Hom(G,V)\to \Hom(G,W)$. Then for each $s\in S$ we have an isomorphism  $\Hom(P_s,V)\to \Hom(P_s,W).$ Therefore using Proposition 2.1.8 of \cite{SchneidersQA}, $V \to W$ is a strict epimorphism. On the other hand, we have $\Hom(P_s,\ker f) = 0$ for every $s\in S$. Therefore $\ker f = 0$ and $f$ is a monomorphism. Any strict epimorphism which is a monomorphism is in fact an isomorphism, so we are done.
\ \hfill $\Box$

In Proposition 2.1.16 (b) of \cite{SchneidersQA} it is shown that:
\begin{prop}\label{prop:FilteringInductive}In elementary quasi-abelian category, filtered inductive limits are exact.
\end{prop}

\subsection{Some Relative Algebraic Geometry}
We will assume that the reader is familiar with the ideas of commutative unital monoids and their modules in a closed symmetric monoidal category with unit object (see \cite{TV}). In this work we will discuss only closed symmetric monoidal categories which also have the structure of quasi-abelian category \cite{SchneidersQA}. We will use notations $\ttComm(\ttC)$ for the category of commutative monoids over $(\ttC, \ootimes, \id_\ttC)$. Also, if $A$ is an object of $\ttComm(\ttC)$ then $\ttMod(A)$ denotes the category of modules over $A$ in the category $\ttC$. We remark that $\ttMod(A)$ is a quasi-abelian closed symmetric monoidal category, whose monoidal functor will be denoted by $\ootimes_A$. So it is meaningful to talk about the bounded above derived category of $\ttMod(A)$, that will be denoted by $D^{\leq 0}(A)$. An important fact about the quasi-abelian structure inherited by $\ttMod(A)$ from $\ttC$ is that the forgetful functor $\ttMod(A) \to \ttC$ commutes with all finite limits and colimits. In particular,  kernels and cokernels can be calculated in $\ttC$ and a morphism in $\ttMod(A)$ is strict if and only if is strict as a morphism of $\ttC$. We will also use the following: 
\begin{lem}\label{lem:2of3discrete} Let $\ttC$ be a quasi-abelian category.
	Suppose that we have an exact triangle in the category $D^{\leq 0}(\ttC)$ 
	\[V_1  \longrightarrow V_2 \longrightarrow V_3.
	\]
	If $V_2$ and $V_3$ are discrete, then $V_1$ is as well.
\end{lem}
\ \hfill $\Box$

The category of affine schemes of $\ttC$ is defined to be the opposite category of $\ttComm(\ttC)$, the functor that associate to a monoid $A$ his opposite is denoted $\spec(A)$ and the category of affine schemes with $\ttAff(\ttC)$.

Let's consider an arbitrary morphism $q:\spec(C) \to \spec(A)$ and the Cartesian diagram 
\begin{equation}\label{basechange}\xymatrix{  \spec(C\ootimes_{A}B) \ar[r]^{q'} \ar[d]_{p'} & \spec(B) \ar[d]^{p} \\
	\spec(C) \ar[r]_{q}& \spec(A). 
}
\end{equation}

There is a natural equivalence 
\begin{equation}p^{*}q_{*} \Longrightarrow q'_{*}p'^{*}
\end{equation}
called base change. 

\begin{defn}A base change of a morphism $p:\spec(B) \to \spec(A)$ is the morphism $p'$ appearing in diagram (\ref{basechange}) for some $q$.
\end{defn}
The following notion appears in \cite{TVe}, \cite{TVe2}, \cite{TVe3}, \cite{TVe4}:
\begin{defn}A morphism $\spec(B)\to \spec(A)$ is a \emph{homotopy monomorphism} if the canonical functor $D^{\leq 0}(B)\longrightarrow D^{\leq 0}(A)$ is fully faithful. In a dual way we say that the correspondent morphism of monoids $A \to B$ is a \emph{homotopy epimorphism}.
\end{defn}

\begin{lem} \label{lem_HomotopyMon}
	Assume that $p:\spec(B)\to \spec(A)$ is a morphism in $\ttAff(\ttC)$ and that the functor $\ttMod(A) \to \ttMod(B)$ given by tensoring with $B$ over $A$ is explicitly left derivable to a functor $D^{\leq 0}(A)\to D^{\leq 0}(B)$. Then $p$ is a homotopy monomorphism if and only if $B\ootimes^{\mathbb{L}}_{A} B\cong B$.
\end{lem}

{\bf Proof.} 
For any object $M$ of $D^{\leq 0}(B)$ we have 
\begin{equation*}
M\ootimes^{\mathbb{L}}_{A}B\cong M\ootimes^{\mathbb{L}}_{B}(B\ootimes^{\mathbb{L}}_{A}B).
\end{equation*}

Hence $\mathbb{L}p^{*}p_{*}\to  id_{D^{\leq 0}(B)}$ is an isomorphism if and only if  we have natural isomorphisms $M\ootimes^{\mathbb{L}}_{A}B\cong M$ for any $M \in D^{\leq 0}(B)$ which happens if and only if $B\ootimes^{\mathbb{L}}_{A}B\cong B$.
\ \hfill $\Box$

The following is an easy consequence of the definitions:

\begin{lem}\label{lem:exchange}For any epimorphism $A\to B$ and any $B$-modules, $M$ and $N$, the natural morphism $M \ootimes_{A}N\to M \ootimes_{B} N$ is an isomorphism. For any homotopy epimorphism $A\to B$ and any  $M, N\in D^{\leq 0}(A)$, the natural morphism $M \ootimes^{\mathbb{L}}_{A}N\to M \ootimes^{\mathbb{L}}_{B} N$ is an isomorphism.
\end{lem}
{\bf Proof}
If $A\to B$ is an epimorphism, then $B\ootimes_{A}B \to B$ is an isomorphism so 
\[M \ootimes_{A} N \cong (M\ootimes_{B} B) \ootimes_{A} (B \ootimes_{B} N) \cong M \ootimes_{B} (B \ootimes_{A} B) \ootimes_{B} N\to M \ootimes_{B} B \ootimes_{B} N \cong M \ootimes_{B}N
\]
is an isomorphism.
If $A\to B$ is a homotopy epimorphism, then $B\ootimes^{\mathbb{L}}_{A}B \to B$ is an isomorphism so 
\[M \ootimes^{\mathbb{L}}_{A} N \cong (M\ootimes^{\mathbb{L}}_{B} B) \ootimes^{\mathbb{L}}_{A} (B \ootimes^{\mathbb{L}}_{B} N) \cong M \ootimes^{\mathbb{L}}_{B} (B \ootimes^{\mathbb{L}}_{A} B) \ootimes^{\mathbb{L}}_{B} N\to M \ootimes^{\mathbb{L}}_{B} B \ootimes^{\mathbb{L}}_{B} N \cong M \ootimes^{\mathbb{L}}_{B}N
\]
is an isomorphism.
\ \hfill $\Box$


\begin{defn}\label{defn:flat}Let $(\ttC, \ootimes, \id_{\ttC})$ be a closed, symmetric monoidal, quasi-abelian category. We call an object $F$ of $\ttC$ \emph{flat} if for any strictly exact sequence 
	\[0 \to E' \to E \to E'' \to 0
	\]
	the resulting sequence 
	\[0 \to E'\ootimes F \to E\ootimes F \to E'' \ootimes F \to 0
	\]
	is strictly exact, \ie if the endofunctor $E \mapsto E \ootimes F$ is an exact functor.
\end{defn}
The following is Proposition 2.1.19 of \cite{SchneidersQA}, which tell us that the construction of Ind-categories behave well for quasi-abelian closed symmetric monoidal categories. 

\begin{prop}
	Let $\mathcal{E}$ be a small quasi-abelian category with enough flat projectives. Assume $\mathcal{E}$ is also a closed symmetric monoidal category, then the same holds for $\ttInd(\mathcal{E})$ in an essentially unique way extending the structure on $\mathcal{E}$. If for every projective object $P$ of $\mathcal{E}$ we have that 
	\[ \mathcal{E}\to \mathcal{E}
	\]
	\[E \mapsto P\ootimes E
	\]
	is exact and that $P \ootimes P'$ is projective for any projective object $P'$ of $\mathcal{E}$, then the same properties hold in $\ttInd(\mathcal{E})$.
\end{prop}

The following is an elaboration of a special case of Remark 1.3.21  of \cite{SchneidersQA}.
\begin{lem}Let $(\ttC, \ootimes, \id_{\ttC})$ be a closed symmetric monoidal quasi-abelian category with enough flat projectives. Then any projective object of $\ttC$ is flat. If in addition there are enough projectives in $\ttC$ then the full additive sub-category of projectives in $\ttC$ is projective in the sense of Definition 1.3.2 of \cite{SchneidersQA} for the functor $W \mapsto W \ootimes V$ for any $V\in \ttC$. This allows us to define the left derived functors of such functors.
\end{lem}
{\bf Proof.} Let $P$ be a projective object of $\ttC$. Consider a strictly exact sequence 
\[0 \to E' \to E \to E'' \to 0.
\]
Choose a strict epimorphism $F \to P$ where $F$ is flat. Since $P$ is projective, we can split it so that $F \cong P \oplus P'$ and hence, because of the assumptions on $F$, we also have that 
\[0 \to E'\ootimes P\to E\ootimes P\to E'' \ootimes P \to 0
\]
is strictly exact. Suppose that we are given a strictly exact sequence 
\[0 \to E' \to E \to E'' \to 0.
\]
where $E$ and $E''$ are projective. Then, since $E''$ is projective $E \cong E' \oplus E''$. Given any strict epimorphism $X \to Y$ we can split the map $\Hom(E, X) \to \Hom(E, Y)$ (surjective since $E$ is projective) into maps $\Hom(E', X) \to \Hom(E', Y)$  and $\Hom(E'', X) \to \Hom(E'', Y)$ which are therefore surjective. Hence $E'$ is projective. Therefore  the full additive sub-category of projectives in $\ttC$ is projective for the functor $W \mapsto W \ootimes V$.
\ \hfill $\Box$
\begin{lem}Let $(\ttC, \ootimes, \id_{\ttC})$ be a closed symmetric monoidal quasi-abelian category with enough flat projectives. The flat objects agree with those objects $F$ for which $F \ootimes^{\mathbb{L}}E \to F \ootimes E$ is an isomorphism for every object $E$. This additive sub-category of flat objects is a projective class for the functor $W \mapsto W \ootimes V$ for any object $V$ in the sense of Definition 1.3.2 of \cite{SchneidersQA}.
\end{lem}
{\bf Proof.} In order to show that the flat objects form a projective class for the functors  $W \mapsto W \ootimes V$ , the only non-obvious part of the definition to check is a two out of three rule for flat objects. Consider a strictly exact sequence 
\[0 \to E' \to E \to E'' \to 0
\]
for which $E$ and $E''$ are flat. We get an exact triangle 
\[E' \ootimes^{\mathbb{L}} V \to E  \ootimes^{\mathbb{L}} V\to E''  \ootimes^{\mathbb{L}} V
\]
and because $E  \ootimes^{\mathbb{L}} V\to E \ootimes V$ and $E''  \ootimes^{\mathbb{L}} V\to E''  \ootimes V$ are isomorphisms we see, using Lemma \ref{lem:2of3discrete}, that  $E'  \ootimes^{\mathbb{L}} V\to E' \ootimes V$ is also an isomorphism. The proofs of the remaining statements are the same as in standard homological algebra.
\hfill $\Box$

\subsection{Topologies}\label{Topologies}
Let $(\ttC, \ootimes, \id_{\ttC})$ be a closed symmetric monoidal elementary quasi-abelian category with enough flat projectives. In \cite{BeKr2} we discuss a (Grothendieck) topology of (homotopy) Zariski open immersions on $\ttComm(\ttsC)^{op}$ where $\ttsC$ is the closed symmetric monoidal model category of simplicial objects in $\ttC$. The model structure is compatible in a natural way with the quasi-abelian structure on $\ttC$. This allows us to use the work of To\"{e}n and Vezzosi from \cite{TVe2} and \cite{TVe3}. The covers in this topology consist of collections $\{\spec (B_i)\to \spec(A)\}_{i\in I}$ where there exists a finite subset $J \subset I$ such that
\begin{itemize}
\item
for each $i\in J$, the morphism $A\to B_i$ is homotopically of finite presentation and the resulting morphisms $D^{\leq 0}(B_i)\to D^{\leq 0}(A)$ is fully faithful 
\item
a morphism in $D^{\leq 0}(A)$  between compact objects is an isomorphism if and only if it becomes an isomorphism in each $D^{\leq 0}(B_j)$ for $j\in J$ after applying the functor $M \mapsto M\ootimes_{A}^{\mathbb{L}} B_{j}$. Such a family is called \emph{conservative}.
\end{itemize}
There are many other, non-equivalent, topologies that can be described in this abstract setting. For instance, the first condition can be replaced by flat epimorphisms, \'{e}tale morphisms, smooth morphisms, etc. The second condition corresponds geometrically to each point in the target being in the image of one of the spaces in the cover, it can also be modified in various ways as well. The requirement of finite covers also should probably be changed in some contexts. We cannot emphasize enough that, even when $A$ and $B_i$ are discrete (so objects of $\ttComm(\ttC)$), the morphisms $A\to B_i$ appearing in the covers defined above in our examples do not exhibit the objects $B_i$ as flat in the category $(\ttMod(A), \ootimes_{A}, A)$ in the sense of Definition \ref{defn:flat}. The algebraic tensor product and the completed tensor product can only be expected to agree for finite $A$-modules and all the important examples of localizations $A\to B$ that we know about need not present $B$ as flat over $A$ for the completed tensor product. This was first pointed out to use by A. Ducros. Our use of the word flat does not agree with that of To\"{e}n and Vezzosi and also does not agree with many points in the literature on analytic geometry where the algebraic tensor product is unfortunately used. We believe that stalk-wise flatness (see for instance \cite{Du}) can also be expressed in our language but we have not done so yet. In this article we mainly focus on the first condition, looking at the particular cases such as $\ttC=\ttInd(\ttBan^{A}_{\mathbb{Z}})$ and $\ttC=\ttBorn_{k}$. In \cite{BeKr} we looked at the case $\ttC = \ttBan^{nA}_{k}$. The conceptual idea is to show that on sub-categories of $\ttComm(\ttC)^{op} \subset \ttComm(\ttsC)^{op}$ the abstractly defined topologies on  $\ttComm(\ttsC)^{op}$ restrict to well known topologies in analytic geometry. This provides an extension of standard analytic geometry to bigger categories where all kinds of limits and colimits of affine schemes exist.
\section{Banach and Ind-Banach modules}

\subsection{Normed and Banach modules}

We begin this section by recalling some basic facts on normed rings and normed modules for which sometime there is a lack a of precise bibliographic reference and also to fix notations, before studying the category of inductive systems of Banach modules.

\begin{defn}
By a \emph{complete normed (or Banach) ring} we mean a commutative ring with identity $R$ equipped with a function, called \emph{absolute value}, $|\cdot|: R \to \Rb_{\ge 0}$ such that
\begin{itemize}
\item $|a| = 0$ if and only if $a=0$;
\item $|a + b| \le |a| + |b|$ for all $a,b \in R$;
\item there is a $C>0$ such that $|a b| \leq C |a||b|$ for all $a,b \in R$;
\item $R$ is a complete metric space with respect to the metric $(a,b) \mapsto |a-b|$.
\end{itemize}
The category of complete normed rings has as objects complete normed rings and as morphisms \emph{bounded homomorphisms} of rings, \ie ring homomorphisms $R \to S$ such that there exists a constant $C > 0$ such that $|\phi(a)|_B \le C |a|_A$ for all $a \in A$. A Banach ring is called \emph{non-Archimedean} if its absolute value satisfies the strong triangle inequality.
 

\end{defn}
\begin{defn}\label{defn:seminrm}
Let $(R, |\cdot|_R)$ be a complete normed ring. A \emph{semi-normed module} over $R$ is an $R$-module $M$ equipped with a function $\|\cdot \|_M: M \to \Rb_{\ge 0}$ such that for any $m,n \in M$ and $a \in R$:
\begin{itemize}
\item $\|0_M\|_M = 0$;
\item $\|m + n\|_M \le \|m\|_M + \|n\|_M$;
\item $\|a m\|_M \le C |a|_R\|m\|_M$ for some constant $C > 0$.
\end{itemize}
A \emph{normed} module is a semi-normed module for which $\|m\|_M = 0$ implies that $m=0_M$ and a \emph{Banach} module is a normed module for which every Cauchy sequence converges. 
\end{defn}
\begin{defn}
A complete normed ring or a semi-normed module over a complete normed ring is called \emph{non-Archimedean} if its semi-norm obeys the strong triangle inequality.
\end{defn}

\begin{defn}Let $(R, |\cdot|_R)$ be a complete normed ring. 
A homomorphism between semi-normed $R$-modules, $f: (M, \|\cdot\|_M) \to (N, \|\cdot\|_N)$ is called \emph{bounded} if there exists a real constant $C > 0$ such that
\[ \|f(m)\|_N \le C \|m\|_M \]
for any $m \in M$.
\end{defn}

From now on we will always suppose that $(R, |\cdot|_R)$ is a Banach ring.

\begin{notation}
We will use the following notation:
\begin{itemize}
\item $\ttSNrm^{A}_R$ the category of semi-normed modules;
\item $\ttNrm^{A}_R$ the category of normed modules;
\item $\ttBan^{A}_R$ the category of Banach modules.
\end{itemize}
where the hom-sets are always the sets of bounded homomorphisms between $R$-modules and we suppress an explicit reference to the norm of $R$.
We use $V \mapsto V^{s}$ to be the left adjoint to the inclusion $\ttNrm^{A}_R \subset \ttSNrm^{A}_{R}$ and $V\mapsto \widehat{V}$ to be the left adjoint to the inclusion $\ttNrm^{A}_R \subset \ttBan^{A}_{R}$. For $V \in \ttSNrm^{A}_{R}$ we often write $\widehat{V}$ instead of $\widehat{V^{s}}$. The ``$A$" in superscript will be motivated later when we will introduce the categories of non-Archimedean modules over $R$.
\end{notation}

\begin{defn}
Two semi-norms $\|\cdot\|, \|\cdot\|'$ on a abelian group $M$ are said to be \emph{equivalent} if there exist two constants $c, C > 0$ such that
\[ c \|\cdot\|' \le \|\cdot\| \le C \|\cdot\|'. \]
\end{defn}

Hence, two semi-norms on a $R$-modules $M$ are equivalent if and only if the identity map $M\to M$ induces an isomorphism $(M, \|\cdot\|) \cong (M, \|\cdot\|')$ in $\ttSNrm^{A}_R$.

\begin{rem} \label{rem:endo}
The reader may have wondered why we did not take $C = 1$ in the Definition \ref{defn:seminrm} of a semi-normed module. In fact if $(M, \|\cdot\|_M)$ is a semi-normed module over $(R, |\cdot|_R)$ we can consider the semi-norm $\|\cdot\|_M': M \to \Rb_{\ge 0}$ defined by
\[ \|m\|_M' = \sup_{a \in R, a \ne 0} \left \{ \frac{\|a m\|_M}{|a|_R} \right \}. \]
This always gives a finite number because
\[
\|m\|_M' \le \sup_{a \in R, a \ne 0} \left \{ \frac{C |a|_R \|m\|_M}{|a|_R} \right \} \le C \|m\|_M  , \]
thus this semi-norm is always well defined.
Then $\|\cdot\|_M'$ is a semi-norm on $M$, equivalent to $\|\cdot\|_M$ (to be shown in Proposition \ref{prop:endo}) for which $\|a m\|_M' \le |a|_R\|m\|_M'$. Notice that a morphism of $R$-modules from $M$ to $N$ is a bounded morphism from $(M, \|\cdot\|_M)$ to $(N, \|\cdot\|_N)$ if and only if the same morphism of $R$-modules defines a bounded morphism $(M, \|\cdot\|_M')$ to $(N, \|\cdot\|_N')$. 
\end{rem}

So, we can always suppose to have a module $(M, \|\cdot\|_M)$ equipped with a semi-norm such that $\|a m\|_M \le |a|_R\|m\|_M$, and we will see in Proposition \ref{prop:endo} that nothing that follows will depend on this choice. 

\begin{lem} 
Let $(M, \|\cdot\|_M)$ be an object of $\ttSNrm^{A}_R$ and $\pi: M \to M/N$ be the canonical quotient morphism, then the residue semi-norm
\[ \|q\|_{M/N} = \inf_{m \in q} \|m\|_M \]
makes $M/N$ into an semi-normed $R$-module, and $\pi$ is a strict epimorphism.

\end{lem}

{\bf Proof.}

The only non-obvious thing to check is that the residue semi-norm is well defined. And to do this is enough to check that 
\[ \|a q\|_{M/N} \le |a|_R \|q\|_{M/N} \]
for any $q \in M/N$ and $a \in R$. Indeed
\[ \|a q\|_{M/N} = \inf_{m \in a q} \|m\|_M = \inf_{m \in q} \|a m\|_M \le 
  \inf_{m \in q} |a|_R \|m\|_M = |a|_R \|q\|_{M/N}. \]
\ \hfill $\Box$

\begin{rem}
Notice that in the case when $R$ is a valued field the conditions $\|a m\|_M = |a|_R\|m\|_M$ and $\|a m\|_M \le |a|_R\|m\|_M$ are equivalent and this explain why one can works with the former definition. If $M$ is equipped with a semi-norm such that $\|a m\|_M = |a|_R\|m\|_M$ holds, then $M$ is often said to be \emph{faithfully semi-normed}. The category of faithfully semi-normed modules is well behaving only over valued fields, otherwise (for example) quotients of faithfully semi-normed modules are not faithfully semi-normed modules.
\end{rem}

From now on we will consider only $\ttBan^{A}_R$, but a similar study can be done also for $\ttSNrm^{A}_R$ and $\ttNrm^{A}_R$.

\begin{prop} 
Let $f: M \to N$ be a morphism in $\ttBan^{A}_R$, then
\begin{itemize}
\item $\ker(f) = f^{-1}(0)$ equipped with the restriction of the norm of $M$;
\item $\coker(f) = N/\ol{f(M)}$ equipped with the residue norm.
\end{itemize}
{\bf Proof.}
Easy verifications.
\ \hfill $\Box$
\end{prop}

\begin{prop} 
The category $\ttBan^{A}_R$ is pre-abelian and the coproduct of a finite family $(M_1, \|\cdot\|_{M_1}), \ldots, (M_n, \|\cdot\|_{M_n})$ of Banach modules is given by $(M_1 \oplus \ldots \oplus M_n, \|\cdot\|_{M_1 \oplus \ldots \oplus M_n} )$ where
\[ \|\cdot\|_{M_1 \oplus \ldots \oplus M_n} = \sum_{i = 1}^n \|\cdot\|_{M_i}. \]
\end{prop}

{\bf Proof.}
It is easy to see that $(M_1 \oplus \ldots \oplus M_n, \|\cdot\|_{M_1 \oplus \ldots \oplus M_n} )$ satisfies the required universal property.
\hfill $\Box$

\begin{rem}
On the direct sum $M_1 \oplus \ldots \oplus M_n$ of the underlying $R$-modules (forgetting the norms) one can also consider the norm
\[ \|(m_1, \ldots, m_n)\|_{\max{}}  = \max_{1 \le i \le n} \{ \|m_i\|_{M_i} \}. \]
This norm defines the product of the $M_i$. However, this norm turns out to be equivalent to $|\cdot|_{M_1 \oplus \ldots \oplus M_n}$ because
\[ \|(m_1, \ldots, m_n)\|_{\max{}} \le \|(m_1, \ldots, m_n)\|_{M_1 \oplus \ldots \oplus M_n} \le n \|(m_1, \ldots, m_n)\|_{\max{}}  \]
for any $(m_1, \ldots, m_n) \in M_1 \oplus \ldots \oplus M_n$. Which shows that the norm on finite products are essentially unique from which descent the pre-abelian nature of $\ttBan^{A}_R$. We will say more about products and coproducts in the end of this section.
\end{rem}

\begin{prop} 
Let $f: M \to N$ be a morphism in $\ttBan^{A}_R$, then
\begin{itemize}
\item $\im(f) \cong \ol{f(N)}$ with the norm induced by $N$;
\item $\coim(f) \cong M/\ker(f)$ with the residue norm;
\item $f$ is a monomorphism if and only if is injective;
\item $f$ is a strict monomorphism if and only if is injective and the norm on $M$ is equivalent to the norm induced by $N$ and $M$ is a closed subset of $N$;
\item $f$ is a strict epimorphism if and only if is surjective and the residue norm of $M/\ker(f)$ is equivalent to the semi-norm of $N$.
\end{itemize}
\end{prop}
{\bf Proof.}
The proofs of these results are similar enough to \cite{BeKr} and \cite{Pr} that we do not repeat them in this context here.

\ \hfill $\Box$

We underline that a strict morphism $f: M \to N$ in $\ttBan^{A}_R$ might be not a strict morphism in $\ttSNrm^{A}_R$, because in $\ttSNrm^{A}_R$ one has that $\im(f) \cong f(N)$. But if one consider an exact sequence
\[ 0 \to M' \to M \to M'' \to 0 \]
in $\ttBan^{A}_R$ then the sequence is strict in $\ttBan^{A}_R$ if and only if is strict in $\ttSNrm^{A}_R$

\begin{prop} \label{prop:quasi_abelian}
The category $\ttBan^{A}_R$ is a quasi-abelian category.
\end{prop}
{\bf Proof.}
We have to show that pullbacks preserve strict epimorphisms and pushouts preserve strict monomorphisms. We want to underline that the proof follows very closely the proof that Schneiders wrote for the case when $R$ is a valued field (see section 3.2.1 of \cite{SchneidersQA}), but is not completely identical to that.  So, let
\[
\begin{tikzpicture}
\matrix(m)[matrix of math nodes,
row sep=2.6em, column sep=2.8em,
text height=1.5ex, text depth=0.25ex]
{M & U  \\
 V & W \\};
\path[->,font=\scriptsize]
(m-1-1) edge node[auto] {$p$} (m-1-2);
\path[->,font=\scriptsize]
(m-1-1) edge node[auto] {$q$} (m-2-1);
\path[->,font=\scriptsize]
(m-2-1) edge node[auto] {$g$}  (m-2-2);
\path[->,font=\scriptsize]
(m-1-2) edge node[auto] {$f$}  (m-2-2);
\end{tikzpicture}
\]
be a Cartesian square where $g$ is a strict epimorphism. We have to show that $p$ is a strict epimorphism. Defining the map
\[ \a: U \oplus V \to W, \ \ \a = (f, - g) \]
then 
\[ M \cong \ker(\a). \]
So, $g$ is surjective and is easy to check that $q$ is surjective and also that if we denote $i: \ker(\a) \to U \oplus V$ the canonical morphism, then
\[ q = \pi_V \circ i,  \ \ p = \pi_U \circ i. \]
Moreover, recall that
\[ \|\cdot\|_{U \oplus V} = \|\cdot \|_U + \|\cdot \|_V \]
and that $\ker(\a)$ is equipped with the restriction of $\|\cdot \|_{U \oplus V}$, so we have to show that
\[ U \cong \frac{M}{\ker(q)}. \]
Clearly
\[ \ker(p) = \{ (x, y) \in \ker(\a) = M | x = 0, g(y) = 0 \} \]
so for $(x,y) \in M$ we have
\[ \inf_{(u, v) \in \ker p} \|(x, y) + (u , v)\|_{U \oplus V} = \inf_{v \in \ker g} \|y + v\|_V + \|x\|_U. \]
And by hypothesis $g$ is a strict epimorphism, so there exists $r > 0$ such that
\[ \inf_{v \in \ker g} \|y + v\|_V \le r \|g(x)\|_W \]
for any $y \in V$. Then, since $f$ is a bounded morphism there exist a constant $r' > 0$ such that
\[ \|g(y)\|_W =  \|f(x)\|_W \le r' \|x\|_U \]
for any $(x, y) \in M$. Hence
\[ \inf_{(u, v) \in \ker p} \|(x, y) + (u , v)\|_{U \oplus V} \le r r' \|x\|_U + \|x\|_U = (1 + r r') \|x\|_U \]
which shows that $p$ is a strict epimorphism.

Now we have to check the dual statement, \ie that pushouts preserve strict monomorpshims. So let
\[
\begin{tikzpicture}
\matrix(m)[matrix of math nodes,
row sep=2.6em, column sep=2.8em,
text height=1.5ex, text depth=0.25ex]
{W & U  \\
 V & M \\};
\path[->,font=\scriptsize]
(m-1-1) edge node[auto] {$g$} (m-1-2);
\path[->,font=\scriptsize]
(m-1-1) edge node[auto] {$f$} (m-2-1);
\path[->,font=\scriptsize]
(m-2-1) edge node[auto] {$p$}  (m-2-2);
\path[->,font=\scriptsize]
(m-1-2) edge node[auto] {$q$}  (m-2-2);
\end{tikzpicture}
\]
be a co-Cartesian square where $g$ is a strict monomorphism. We have to check that $p$ is a strict monomorphism. First, let's show that the norm on $V$ is equivalent to the norm induced by $M$.

We can put
\[ M = \coker(\a: W \to U \oplus V) \cong \frac{U \oplus V}{\overline{\a(W)}} \]
equipped with the residue norm. Let's denote 
\[ \pi: U \oplus V \to \coker(\a) \]
the canonical map, and we remark that
\[ p = \pi \circ i_V, q = \pi \circ i_U.  \]
It is easy to see that $p$ is injective, let's check that it is strict. Since $f$ is bounded and $g$ is strict we can find positive constants $r, r', r'', r''' > 0$ such that
\begin{eqnarray*}
\|y\|_V & \le & \|y + f(x) - f(x)\|_V  \\
      & \le & \|y + f(x)\|_V + r \|f(x)\|_V \\
      & \le & \|y + f(x)\|_V + r' \|x\|_W  \\
      & \le & \|y + f(x)\|_V + r'' \|g(x)\|_U \\ 
      & \le & r''' \| (g(x), y + f(x)) \|_{U \oplus V}
\end{eqnarray*}
for all $y \in V$ and $x \in W$. Hence for any $y \in V$
\begin{eqnarray*}
\|y\|_V & \le & \inf_{x \in W} r''' \| (g(x), y + f(x)) \|_{U \oplus V}  \\
      & \le & \inf_{(u,v) \in \a(W)} r''' \| (u, y + v) \|_{U \oplus V} \\
      & \le & \inf_{(u,v) \in \a(W)} r''' \| q((0, y)) + (u, v) \|_{U \oplus V}  \\
      & \le & \inf_{(u,v) \in \a(W)} r''' \| p(y) + (u, v) \|_{U \oplus V} 
\end{eqnarray*}
which shows that $p$ is a strict monomorphism.

Now, it remain to show that $p(V)$ is closed in $M$. By looking at the strictly coexact sequence in $\ttBan^{A}_R$
\[ W \stackrel{\binom{g}{-f}}{\to} U \oplus V \stackrel{\pi}{\to} M \to 0  \]
one can deduce that this sequence is also strictly coexact as a sequence of $\ttSNrm^{A}_R$. Let $\{x_n\}$ be a sequence of elements of $V$ such that $\lim_{n \to \infty} p(x_n) = y \in M$. Setting $y = \pi(u, v)$, we have that
\[ \pi((0, x_n) - (u, v)) = \pi((-u, x_n - v)) \to 0 \] 
as $n \to \infty$. And since the above sequence is strictly coexact there is a sequence $\{ w_n \}$ of $W$ such that 
\[ (-u - g(w_n), x_n - v + f(w_n)) \to 0 \]
in $U \oplus V$. Then $g(w_n) \to -u$ in $U$, hence $w_n \to w$ so that $g(w) = -u$. Finally 
\[ x_n \to x = v - f(w_n) = v - f(w) \]
in $V$ which shows that
\[ p(x) = \pi((0, x)) = \pi((u, v) + (-g(w), f(w))) = y \]
which concludes the proof.

\ \hfill $\Box$

$\ttBan^{A}_R$ is equipped with an internal hom-functor defined as follows: given $M,N \in \ttBan^{A}_R$, then the $R$-module $ \Hom(M, N)$ is equipped with the semi-norm defined by
\[ \|f\|_\sup = \sup_{m \in M, \|m\|_M \ne 0} \frac{\|f(m)\|_N}{\|m\|_M}.  \]  In fact, it is easy to check that it is non-degenerate and complete and we consider it as an object $\uHom_{\ttBan^{A}_R}(M, N)$ of $\ttBan^{A}_R$. 
Let $M, N \in \ttBan^{A}_R$ then the underlying $R$-module of $M \otimes_R N$ can be quipped with the semi-norm
\begin{equation}\label{eqn:DefnSnrm} \|x\|_{M \otimes_{R} N} = \inf \left \{ \sum_{i\in I} \|m_i\|_M \|n_i\|_N \ \ | \ \ x = \sum_{i\in I} m_i \otimes n_i \ \ , \ \ |I|<\infty \right \}. \end{equation}
This is called their \emph{projective (Archimedean) tensor product}. Notice in general, $(M \otimes_R N, \|\cdot\|_{M \otimes_{R} N})$ is not an object of $\ttBan^{A}_R$, 
\begin{defn}\label{def:Archimedean_tensor_product}
The monoidal structure we use on $\ttBan^{A}_{R}$ is the so called \emph{complete projective (Archimedean) tensor product} defined by
\[ M \wotimes_R N = \widehat{M \otimes_R N}\]
where we preform the (seperated) completion of $M\otimes_{R}N$ equiped with the semi-norm from (\ref{eqn:DefnSnrm}). 
\end{defn}

The following Proposition, which shows that the projective tensor product is well defined, is a well know result, but can be difficult to find in literature.

\begin{prop}\label{prop:WD}
Let $M, N \in \ttBan^{A}_R$ then for any $\la \in R$ and any $x \in M \wotimes_R N$
\[ \|\la x\|_{M \wotimes_{R} N} \le |\la|_R \|x\|_{M \wotimes_{R} N}. \]
\end{prop}
{\bf Proof.}
Consider an expression $x =\sum_{i\in I} m_{i} \otimes n_{i}$ where $|I |<\infty$.
We have 
\begin{equation}
\begin{split} 
 \|\la x\|_{M \otimes_{R} N} & = \inf \left \{ \sum_{i\in I} \|m_i\|_M \|n_i\|_N \ \ |  \ \ \la x = \sum_{i\in I} m_i \otimes n_i \right \} \\ & \le 
 \inf \left \{ |\la| \sum_{i\in I} \|m_i\|_M \|n_i\|_N \ \ | \ \ x = \sum_{i\in I} m_i \otimes n_i \right \} = |\la|_R \|x\|_{M \otimes_{R} N}
\end{split}
\end{equation}

because if $x = \sum_{i\in I} m_i \otimes n_i$ then 
\[ \la \sum_{i\in I} m_i \otimes n_i = \sum_{i\in I} (\la m_i) \otimes n_i = \sum_{i\in I} m_i \otimes (\la n_i) = \la x 
\]
hence the first infimum runs over a bigger set of representations of $\la x$, which explains the $\le$ sign. The required inequality for the completed tensor product follows by passing to Cauchy sequences.
\ \hfill $\Box$

\begin{prop}\label{prop:BanAcsm} 
The category $\ttBan^{A}_R$ equipped with the completed projective tensor product defined in \ref{def:Archimedean_tensor_product} and internal hom described above is a closed symmetric monoidal category.
\end{prop}

{\bf Proof.}

We only check that there is an isomorphism 
\[ \Hom(M \wotimes_{R} N, L) \cong \Hom(M, \uHom(N, L)) \]
for any $M, N, L \in \ttBan_R$, which is the less trivial fact. By the universal property of the algebraic tensor product given any bilinear map $\phi: M \times N \to L$ we can find a linear map $\psi: M \otimes_{R} N \to L$ such that the diagram
\[
\begin{tikzpicture}
\matrix(m)[matrix of math nodes,
row sep=2.6em, column sep=2.8em,
text height=1.5ex, text depth=0.25ex]
{ M \times N &  & M \otimes_{R} N  \\
  & L \\};
\path[->,font=\scriptsize]
(m-1-1) edge node[auto] {$\tau$} (m-1-3);
\path[->,font=\scriptsize]
(m-1-1) edge node[auto] {$\phi$} (m-2-2);
\path[->,font=\scriptsize]
(m-1-3) edge node[auto] {$\psi$}  (m-2-2);
\end{tikzpicture}
\]
is commutative, where $\tau$ is the canonical (bounded) map. Then if $\phi$ is bounded, \ie let's suppose that $\|\phi(x)\|_L \le C \|x\|_{M \times N}$ for some constant $C > 0$, we want to show that also $\psi$ is bounded. Let $f \in M \otimes_{R} N$ and $f = \sum_{i\in I} \a_i \otimes \be_i$ a representation of $f$, then
\[ \psi(\sum_{i\in I} \a_i \otimes \be_i) = \sum_{i\in I} \phi(\a_i, \be_i) \]
hence
\[ \|\psi(f)\|_L = \|\sum_{i\in I} \phi(\a_i, \be_i)\|_L \le \sum_{i\in I} \|\phi(\a_i, \be_i)\|_L \le  C \sum_{i\in I} \|(\a_i, \be_i)\|_{M \times N} \]
and this is true for any possible representation of $f$, hence 
\[ \|\psi(f)\|_L \le \|f\|_{M \otimes_{R} N}. \]
This result extends to the completion of the tensor product by functoriality of the completion.

\ \hfill $\Box$

\begin{prop}\label{prop:CopAgree} 
The monoidal product of $(R, |\cdot|_R)$-algebras, thought of as objects of $\ttBan^{A}_R$ coincides with the underlying object in $\ttBan^{A}_R$ of their coproduct in the category $\ttComm(\ttBan^{A}_R)$. 
\end{prop}
{\bf Proof.}

Let $(R_1, |\cdot|_{R_1})$ and $(R_2, |\cdot|_{R_2})$ be two Banach $R$-algebras. It is enough to check that $(R_1 \widehat{\otimes}_{R} R_2, |\cdot|_{R_1  \widehat{\otimes}_{R}  R_2})$ is a normed ring and the only non-obvious thing to check is the existence of a constant $C > 0$ such that
\[ |x y|_{R_1  \widehat{\otimes}_{R}  R_2} \le C |x|_{R_1  \widehat{\otimes}_{R}  R_2} |y|_{R_1  \widehat{\otimes}_{R}  R_2} \]
for any $x, y \in R_1  \widehat{\otimes}_{R}  R_2$. It is enough to show this for finite representations $x = \sum_{i\in I} \a_i \otimes \be_i$ and $y = \sum_{j\in J} \a_j' \otimes \be_j'$ for $|I|<\infty$ and $|J|<\infty$. Now
\begin{equation}
\begin{split} |x y|_{R_1  \widehat{\otimes}_{R}  R_2} & = \inf \{ \sum_{k\in K} |\gamma_k||\delta_k| \ \ | \ \ \sum_{k\in K}  \gamma_k \otimes \delta_k = x y, |K|<\infty \} \\ &\le \sum_{i\in I,j\in J} |\a_i \a_j'||\be_i \be_j'| \le  C_1 C_2 \sum_{i\in I,j\in J} |a_i||a_j'||\be_i||\be_j'|, 
\end{split}
\end{equation}
Therefore
\begin{equation}
\begin{split}
 |x y|_{R_1  \widehat{\otimes}_{R}  R_2} &\le C_1 \inf \{ \sum \sum_{k\in K} |\a_k||\be_k| \ \ | \ \ \sum_{k\in K}  \a_k \otimes \be_k = x, |K|<\infty  \} \\ & \le C_2 \inf \{ \sum_{k\in K}  |\a_k'||\be_k'| | \sum \a_k' \otimes \be_k' = y, |K|<\infty \} \\ &
 = C_1 C_2 |x|_{R_1  \widehat{\otimes}_{R}  R_2} |y|_{R_1  \widehat{\otimes}_{R}  R_2}. 
 \end{split}
 \end{equation}

\ \hfill $\Box$

If the norm of $(R, |\cdot|)$ is non-Archimedean it is standard to consider a different category of Banach modules and rings over $R$ by restricting to the full sub-category 
\[ \ttBan_R^{nA} \rhook \ttBan^{A}_R \] 
of Banach modules whose norm is non-Archimedean. In the following we will use both categories and the notation $\ttBan_R^A$ want to emphasize that we are considering all Banach modules over $R$, allowing also non-ultrametric modules over non-Archimedean base rings. This latter kind of semi-normed modules have almost always been thought to be pathological and excluded a priori from discussions, but this may be not so true.

\begin{prop}Let $R$ be a non-Archimedean Banach ring. The category $\ttBan_R^{nA}$ is quasi-abelian with enough flat projectives. The category $\ttBan_R^{nA}$ is equipped with a structure of closed symmetric monoidal category analogous to the one explained for $\ttBan_R^{A}$.
\end{prop}
{\bf Proof.}
The  monoidal structure is given by the complete \emph{non-Archimedean projective tensor product}, defined as the separated completion of $M \otimes_{R} N$ with respect to
\[ \|x\|_{M \otimes_{R} N} = \inf \left \{ \max_{i \in I} \|m_i\|_M \|n_i\|_N \ \ | \ \ x = \sum_{i \in I} m_i \otimes n_i , \ \ |I|<\infty \right \}. \]
For the quasi-abelian structure, a simple adaptation of arguments of Proposition \ref{prop:quasi_abelian} works, we omit the details (see also \cite{BeKr} in the case that $R$ is a complete non-Archimedean valuation field).
\ \hfill $\Box$ 

In \cite{BGR} one can find a detailed account of the properties of this category.
\begin{obs} 
The inclusion functor $\ttBan_R^{nA} \rhook \ttBan_R^A$ preserves the internal hom-functors but does not preserve the monoidal structures.
\end{obs}
The categories $\ttBan^{A}_{R}$ and $\ttBan^{nA}_{R}$ do not have infinite limits and colimits and that makes it difficult to construct interesting objects in them. For that, we need to introduce two more closed symmetric monoidal categories.

\begin{defn}\label{def:contracting_cat}
We use $\ttBan_R^{A,\le 1}$ to denote the category of Banach modules over $(R, |\cdot|_R)$ with hom-sets defined by non-expanding morphisms, also called contractions.  When $R$ is non-Archimedean then we use $\ttBan_R^{nA,\le 1}$ to denote the category of  non-Archimedean Banach modules over $(R, |\cdot|_R)$ with hom-sets defined by non-expanding morphisms. 
\end{defn}

\begin{prop}  \label{prop:contracting_limits}
The categories $\ttBan_R^{A,\le 1}$ and $\ttBan_R^{nA,\le 1}$ have all limits and colimits, and moreover
\begin{itemize}
\item the inclusion functor $\ttBan_R^{nA,\le 1} \rhook \ttBan_R^{A,\le 1}$ commutes with limits and finite colimits;
\item the inclusion functor $\ttBan_R^{nA,\le 1} \rhook \ttBan_R^{nA}$ commutes with all finite limits and finite colimits;
\item the inclusion functor $\ttBan_R^{A,\le 1} \rhook \ttBan_R^A$ commutes with finite limits and finite colimits.
\end{itemize}
\end{prop}
{\bf Proof.}
We need only to check that $\ttBan_R^{A,\le 1}$ and $\ttBan_R^{nA,\le 1}$ have all products and coproducts to verify the claim on the existence of all limits. Products are easy to describe, with a uniform description in the Archimedean and non-Archimedean case. For objects $M_i$ indexed by a set $I$ we have
\[\prod_{i}{}^{\leq 1}M_{i} =\{
(m_i)_{i \in I} \in \times_{i \in I}M_{i} \ \ | \ \ \sup_{i \in I} \|m_{i}\| < \infty \}  
\] 
equipped with the norm 
\[ \|(m_i)_{i \in I} \| =\sup_{i \in I} \|m_{i}\|. \]

It is easy to see that if all $M_i$ are Banach, their contracting direct product is a Banach module. 

To describe coproducts, we have to discuss separately the Archimedean and the non-Archimedean cases. In $\ttBan_R^{A,\le 1}$ coproducts looks like the completion of
\[ \left ( \bigoplus_{i \in I} M_i, \|(m_i)_{i \in I} \| = \sum_{i \in I} \|m_i\| \right).
\]

Instead in $\ttBan_R^{nA, \le 1}$ coproducts are given by the completion of
\[ \left ( \bigoplus_{i \in I} M_i, \| (m_i)_{i \in I} \| = \sup_{i \in I} \|v_{i} \| \right ). \]

We will denote this coproducts with
\[ \coprod_{i \in I}{}^{A,\le 1} M_i   \ \ \ \ \text{or} \ \ \ \   \coprod_{i \in I}{}^{nA,\le 1} M_i  . \]

The commutation statements of inclusion functors are clear by the description of limits and colimits.

\ \hfill $\Box$

We end this discussion of Banach modules with some result that we will need to show that the inclusion functors $\ttBan_R^{nA,\le 1} \rhook \ttBan_R^{A,\le 1}$ and $\ttInd(\ttBan_R^{nA}) \rhook \ttInd(\ttBan_R^A)$ have adjoints when the base ring is non-Archimedean.

\begin{lem}\label{lem:boundedNorm}
        Let $\{(M_i, \| \cdot \|_i)\}_{i \in I}$ be a collection of Banach modules.        
    A morphism 
    \[ \phi: \coprod_{i \in I}{}^{\le 1} M_i \to N \] 
    satisfies $\|\phi\|\leq r$ if and only if $\|\phi_i\|\leq r$ for all $i\in I$. This means that for any object $N$ we have that for $r>0$
        \[\Hom(\coprod_{i \in I}{}^{\le 1} M_{i}, N)^{\leq r} = \prod_{i \in I} \Hom(M_i, N)^{\leq r}
        \]
        and similarly 
        \[\Hom(\coker [M_{1} \to M_2], N)^{\leq r} = \Hom(M_2, N)^{\leq r} \cap \ker[ \Hom(M_2, N) \to  \Hom(M_1, N)].
        \]
The statement is true in $\ttBan_R^A$ and $\ttBan_R^{nA}$.
\end{lem}

{\bf Proof.}
Notice that for each $i\in I$, $M_i$ has an obvious morphism into $\coprod_{i \in I}{}^{\le 1} M_i$ by taking the identity in component $i$ and zero elsewhere. A morphism $\phi$ from $\coprod_{i \in I}{}^{\le 1} M_i$ to an object $N \in \ttBan_{R}$, with $\|\phi\| \le r$, is uniquely determined by its components $\phi_{i}$ representing its restriction to $M_i$. And the statement of the lemma can be deduced easily noticing that $M_i \to \coprod_{i \in I}{}^{\le 1} M_i$ is an isometry onto its image. 
\hfill $\Box$

\begin{defn}
        If $R$ is a non-Archimedean Banach ring we define for any $M \in \ttBan^{nA}_{R}$
        \[
        P^{nA}(M)=\{(c_{m})_{m \in M^{\times}} \ \ | \ \ c_{m} \in R, \lim_{m \in M^{\times}} \|c_{m}m\|=0 \}
        \]
        with norm $\|c\| = \sup_{m \in M^{\times}} \|c_{m} m\|$
        and 
        \[\kappa^{nA}_{M}: P^{nA}(M)\to M
        \]
        by 
        \[\kappa^{nA}_{M}(c) = \sum_{v \in M^{\times}} c_{v}v .
        \]
        
        For any Banach ring and $M \in \ttBan^{A}_{R}$ we define
        \[
        P^{A}(M)=\{(c_{m})_{m \in M^{\times}} \ \ | \ \ c_{m} \in R, \sum_{m \in M^{\times}} \|c_{m}m\|<\infty \}
        \]
        with norm $\|c\| = \sum_{m \in M^{\times}} \|c_{m} m\|$. Define
        \[\kappa^{A}_{M}: P^{A}(M)\to M
        \]
        by 
        \[\kappa^{A}_{M}(c) = \sum_{m \in M^{\times}} c_{m}m.
        \]
\end{defn}

\begin{rem} The previous definitions can be restated in term of coproducts the respective non-expanding categories: we have 
        \[P^{nA}(M)= \coprod_{m \in M^{\times}}{}^{nA,\le 1} R_{\|m\|} \in \ttBan^{nA}_{R}
        \] and \[P^{A}(M)= \coprod_{m \in M^{\times}}{}^{A,\le 1} R_{\|m\|} \in \ttBan^{A}_{R},\]
        where $R_{\|m\|}$ is the Banach $R$-module obtained by re-scaling the norm of $R$ by a factor of $\|m\|$.
\end{rem}

\begin{lem}
Let $M \in \ttBan^{A}_{R}$, then $P^{A}(M)$ is a projective, flat object in $\ttBan^{A}_{R}$. The same holds for non-Archimedean base rings and $M \in \ttBan^{nA}_{R}$ with $P^{nA}(M)$.
\end{lem}
{\bf Proof.}
The flatness follows from the fact that the two non-expanding categories $\ttBan^{A, \leq 1}_{R}$ and $\ttBan^{nA, \leq 1}_{R}$ are closed symmetric monoidal categories with the same monoidal structure and internal hom spaces as $\ttBan^{A}_{R}$ and $\ttBan^{nA}_{R}$ respectively, and that both $P^{A}(M)$ and $P^{nA}(M)$ are coproducts of flat objects and hence flat. The projectivity of these objects is proven exactly along the lines of the analogous facts when $R$ is a complete valuation field (see \cite{BeKr}).

\ \hfill $\Box$

\begin{lem}Suppose that $M \in \ttBan^{nA}_{R}$.  Then the canonical morphism $\kappa^{nA}_{M}:P^{nA}(M)\to M$ is a strict epimorphism in $\ttBan^{nA}_{R}$. Similarly, if $M \in \ttBan^{A}_{R}$, then the canonical morphism $\kappa^{A}_{M}:P^{A}(M)\to M$ is a strict epimorphism $\ttBan^{A}_{R}$. Both morphisms have norm less than or equal to $1$. Therefore, the categories $\ttBan_R^A$ and $\ttBan_R^{nA}$ have enough flat projectives.
\end{lem}
{\bf Proof.}

The proof is similar enough to the proof in the case where $R$ is a complete valuation field that we do not repeat it here, see \cite{BeKr}.
\ \hfill $\Box$

\begin{prop} \label{prop:endo}
	The association $(M, \|\cdot\|_M) \mapsto (M, \|\cdot\|_M')$, defined in Remark \ref{rem:endo}, is  an endofunctor of  $\ttBan^{A}_R$ which is an equivalence of categories. It preserves internal the hom functor, and the symmetric monoidal structure. The same holds for the non-Archimedean categories when $R$ is non-Archimedean.
\end{prop}

{\bf Proof.}

In Remark \ref{rem:endo} we saw that the identity map $(M, \|\cdot\|_M) \to (M, \|\cdot\|_M')$ is bounded because $\|\cdot\|_M' \le C \|\cdot\|_M$, where $C$ is the constant appearing in the definition of semi-normed $R$-module. Then, to see that the identity map is bounded also in the other direction is enough to see that the inequality
\[ \frac{\|m\|_M}{|1|_R} \le \|m\|_M' \]
holds. Thus the essential image of the functor $(M, \|\cdot\|_M) \mapsto (M, \|\cdot\|_M')$ is $\ttBan^{A}_R$, showing that this functor is an equivalence. Then, is easy to see that the isomorphism class of the projective tensor product (and therefore its completion) and of the hom sets are invariant under this functor. Explicitly, for $x\in M\otimes_{R}N$,
\begin{equation}
\begin{split}
 \|x\|_{M \otimes_{R} N}' & = \inf \{ \sum \|a_i\|'_M\|b_i\|'_N \ \  | \ \ \sum a_i \otimes b_i = x \} \\
 & \le \inf \{ \sum C_M \|a_i\|_M C_N\|b_i\|_N \ \ | \ \ \sum a_i \otimes b_i = x \} \\ & = 
 C_M C_N \inf \{ \sum \|a_i\|_M \|b_i\|_N \ \ | \ \ \sum a_i \otimes b_i = x \} = C_N C_M \|x\|_{M \otimes N} 
\end{split}
\end{equation}
and 
\[ \|x\|_{M \otimes N}' \ge |1_R|_R^2 \inf \{ \sum \|a_i\|_M \|b_i\|_N | \sum a_i \otimes b_i = x \} = |1_R|_R^2 \|x\|_{M \otimes N}. \]
Then, the sup norm on hom-sets doesn't change because for any $\phi \in \Hom_{\ttBan_R}(M, N)$ corresponding to $\phi' \in \Hom_{\ttBan_R}(M', N')$ one can easily find the estimate
\begin{equation}
  \frac{\| \phi \|}{C_M |1_R|_R} \le \|\phi'\|\ \le C_N |1_R|_R \| \phi \|.
\end{equation}
\ \hfill $\Box$

\subsection{Ind-Banach modules}

In this section we want to study the categories $\ttInd(\ttBan_R^A)$ and $\ttInd(\ttBan_R^{nA})$ and the natural structures they are equipped with. 

\begin{defn}\label{def:KSexact}
        (Definition 3.3.1 of \cite{KS}) Let $F:C \to C'$ be a functor. 
        We say that $F$ is KS right exact if the category $C_{U}$ is filtrant for any $U \in C'$. We say that $F$ is KS left exact is $F^{op}:C^{op} \to C'^{op}$ is KS right exact or equivalently if $C^{U}$ is cofiltrant for any $U \in C'$. We say that $F$ is KS exact if it is both KS right and KS left exact.\end{defn}
See \cite{KS} for the precise definitions of the over and under categories $C_U$ and $C^{U}$ and what it means to be filtrant or cofiltrant.
\begin{lem} \label{lem:exact_direct_limit}
The categories $\ttInd(\ttBan_R^A)$ and $\ttInd(\ttBan_R^{nA})$ are complete and cocomplete elementary quasi-abelian with enough flat projectives and in these categories small filtrant colimits are KS exact in the sense of Definition \ref{def:KSexact}. 
\end{lem}

{\bf Proof.}
Let $\mathcal{E}= ``\colim"_{i \in I} E_{i}$ be an element of $\ttInd(\ttBan^{A}_{R}).$ Define 
\begin{equation}\label{eqn:standardArchIndProj}P^{A}(\mathcal{E}) = \coprod_{i \in I} P^{A}(E_{i})
\end{equation}
where the coproduct is taken in $\ttInd(\ttBan^{A}_{R})$, and 
\begin{equation}\label{eqn:standardNonArchIndProj}P^{nA}(\mathcal{E}) = \coprod_{i \in I} P^{nA}(E_{i})
\end{equation}
where the coproduct is taken in $\ttInd(\ttBan^{nA}_{R})$. Using Lemma \ref{lem:coproducts} we can be more explicit about equations (\ref{eqn:standardArchIndProj}) and (\ref{eqn:standardNonArchIndProj}) if needed. The canonical morphisms \[\kappa^{A}_{\mathcal{E}}:P^{A}(\mathcal{E})\to \coprod_{i \in I}{}^{A} E_{i} \to \mathcal{E}\] and  \[\kappa^{nA}_{\mathcal{E}}:P^{nA}(\mathcal{E})\to \coprod_{i \in I}{}^{nA} E_{i} \to \mathcal{E}\] are then strict epimorphisms. Notice that $P^{A}(\mathcal{E})$ is flat in $\ttInd(\ttBan_R^A)$ and that $P^{nA}(\mathcal{E})$ is flat in $\ttInd(\ttBan_R^{nA})$. \hfill $\Box$

We conclude this section with some concreteness results for the category $\ttInd(\ttBan_R)$. 
\begin{defn} \label{defn:essential_monomorphic}
Recall the following notions:
\begin{itemize}
\item a category $\ttC$ equipped with a functor $U: \ttC \to \ttSets$ is called \emph{concrete} if $U$ is a faithful functor;
\item $\ttC$ is called \emph{concretizable} if such a functor exists;
\item an object $``\limind"_{i \in I} X_i \in \ttInd(\ttC)$ is called \emph{monomorphic} if all the morphisms of the system are monomorphisms;
\item an object $X \in \ttInd(\ttC)$ is called \emph{essentially monomorphic} if it is isomorphic to a monomorphic object.
\end{itemize}
\end{defn}

For any concrete category $(\ttC, U_\ttC)$, we define the functor $U_{\ttInd(\ttC)}: \ttInd(\tC) \to \ttSets$ by
\[ U_{\ttInd(\ttC)}(``\limind_{i \in I}" M_i) = \limind_{i \in I} U_{\ttC}(M_i). \]

\begin{prop}	
	Let $(\ttC, U_\ttC)$ be a concrete category such that monomorphisms of $\ttC$ are precisely morphism for which $U_\ttC(f)$ is injective. Then, the sub-category of essentially monomorphic objects of $\ttInd(\ttC)$ is a concrete category with respect to the functor $U_{\ttInd(\ttC)}$.
\end{prop}
{\bf Proof.}
Let $E, F \in \ttInd(\ttC)$ and consider two different morphisms $f, g: E \to F$, work which we suppose to have a common re-indexing, as give by Proposition \ref{prop:reindexing}. This means that there exists at least one $i$ for which $f_i \ne g_i$, and so an element $x_i \in E_i$ such that $f_i(x_i) \ne g_i(x_i)$. But then, since $f = \limind f_i$ and $g = \limind g_i$ we get that for $U_\ttC(f)(y) \ne U_\ttC(g)(y)$, where $y$ is the image of $x_i$ in $\limind U_\ttC(E_i)$, because monomorphism in $\ttC$ precisely corresponds to injective maps and therefore $U_{\ttInd(\ttC)}(``\limind" E_i) = \bigcup U_\ttC(E_i)$ and $U_{\ttInd(\ttC)}(``\limind" F_i) = \bigcup U_\ttC(F_i)$.
\ \hfill $\Box$

In particular we get the following corollary

\begin{cor}  \label{prop:essential_concrete}
	The full sub-category of essentially monomorphic objects of $\ttInd(\ttBan^{A}_R)$ (or $\ttInd(\ttBan^{nA}_R)$ when $R$ is non-Archimedean) is a concrete category with respect to the restriction of the functor $U$.
\end{cor}

\begin{rem} \label{rem:non_concrete_ban}
The category $\ttInd(\ttBan^{A}_R)$ is not a concrete category with respect to the functor $U$. It is enough to show that there two different morphisms $f, g: ``\limind" E_i \to ``\limind" F_i$ such that $U(f) = U(g)$ for a pair of $\ttInd(\ttBan^{A}_R)$.  Let's consider the following non-essentially monomorphic object: let $E$ be a free Banach module of infinite dimension and let's consider a sequence $\{ e_n \}_{n \in \Nb}$ of linearly independent elements. So, for any $n$ the sub-module  $C_n = \lt e_1, \ldots, e_n \gt$ is closed in $E$ and hence $E_n = E / C_n$ is a Banach module. For any $n \le m$ there is a canonical morphism $E_n \to E_m$ and it can be shown that the inductive system $``\limind" E_n$ is not essentially monomorphic, see \cite{PrSc} Remark 3.6. Then, is clear that $U_{\ttInd(\ttC)}(``\limind" E_n) = E / \lt e_1, \ldots, e_n, e_{n+1}, \ldots \gt$.  Let's define a morphism $f: ``\limind" E_n \to ``\limind" E_n$ in the following way. Let's pick a basis $\{ x_i \}_{i \in I}$ which completes $\{e_n \}_{n \in \Nb}$, hence $\{ x_i \}_{i \in I} \cup \{e_n \}_{n \in \Nb}$ is a basis for $E$. For any $n$, $f_n: E_n \to E_n$ is defined to be the morphism which are defined by their action on the basis: $f_n(x_i) = x_i$ and $f_n(e_m) = e_{m + 1}$, for all $m > n$. 

This morphism is not the identity because $\coker(f) \cong R$, the one dimensional free module, while $U_{\ttInd(\ttC)}(f)$ is the identity map on $U_{\ttInd(\ttC)}(``\limind" E_n)$, proving that $\ttInd(\ttBan^{A}_R)$ is not concrete. The same conclusion holds for $\ttInd(\ttBan^{nA}_R)$ if $R$ is non-Archimedean..
\end{rem}

Then, the following Proposition tells that sub-objects of ''concrete" objects are ''concrete".

\begin{prop} \label{prop:essentially_subobject}
Let $E \in \ttInd(\ttBan^{A}_R)$ be an essentially monomorphic object and $E \to F$ be a monomorphism in $\ttInd(\ttBan^{A}_R)$, then $E$ is essentially monomorphic. The same holds in $\ttInd(\ttBan^{nA}_R)$ if $R$ is non-Archimedean.
\end{prop}
{\bf Proof.}

Let $f: E \to F$ be a monomorphism to an essentially monomorphic object in $\ttInd(\ttBan^{A}_R)$. By Proposition \ref{prop:strict_morphism_ind} we know that there exists a re-indexing of $f$: a small filtered category and a natural transformation $f_i:E_i \to F_i$ for $i \in I$ such that $f$ is given by the corresponding morphism
\[ F \cong ``\limind"_{i\in I} F_i \stackrel{\limind f_i}\to ``\limind"_{i\in I} E_i \cong E \]
and all $f_i$ are monomorphisms. Hence for any $i < j$ there is a commutative diagram
\[
\begin{tikzpicture}
\matrix(m)[matrix of math nodes,
row sep=2.6em, column sep=2.8em,
text height=1.5ex, text depth=0.25ex]
{ F_i & E_i  \\
  F_j & E_j \\};
\path[->,font=\scriptsize]
(m-1-1) edge node[auto] {$f_i$} (m-1-2);
\path[->,font=\scriptsize]
(m-1-1) edge node[auto] {$\phi_{i, j}^F$} (m-2-1);
\path[->,font=\scriptsize]
(m-1-2) edge node[auto] {$\phi_{i,j}^E$}  (m-2-2);
\path[->,font=\scriptsize]
(m-2-1) edge node[auto] {$f_j$}  (m-2-2);
\end{tikzpicture}
\]
where $\phi_{i,j}^E$ is a monomorphism. Hence $f_j \circ \phi_{i,j}^F = \phi_{i,j}^E \circ f_i$ is a monomorphism which implies that $\phi_{i,j}^F$ is a monomorphism, therefore $F$ is isomorphic to a monomorphic object.

\ \hfill $\Box$

\begin{obs}The inclusions $\ttBan^{A}_{R}\to \ttInd(\ttBan^{A}_{R})$ and $\ttBan^{nA}_{R} \to \ttInd(\ttBan^{nA}_{R})$ (when $R$ is non-Archimedean) are fully faithful, respect strict morphisms, monomorphisms and epimorphisms, finite limits and colimits and the internal Hom and monoidal structures. The image of any projective (respectively flat) object under any of these functors is projective (respectively flat). In particular, a morphism in $\ttComm(\ttBan^{A}_{R})$ (respectively $\ttComm(\ttBan^{nA}_{R})$) is a homotopy epimorphism if and only if it is in $\ttComm(\ttInd(\ttBan^{A}_{R}))$  (respectively $\ttComm(\ttInd(\ttBan^{nA}_{R}))$).
\end{obs}
\subsection{Complete bornological vector spaces}

In this section we recall the theory of bornological vector spaces and bornological algebras over a non-trivially valued field $k$, and relate it to the theory of $\ttInd(\ttBan_k)$-spaces and algebras. Our aim is to recall results from \cite{Bam}, where the bornological language is used, and recast them in the language of $\ttInd(\ttBan_k)$-spaces. The reason of this choice, will become clear in final sections of this article where we need to use inductive systems to be able to extend the theory over Banach rings. 

We will use in this subsection the convention that when $k$ is a non-Archimedean valued field, then $\ttBan_k = \ttBan_k^{nA}$ as usual in non-Archimedean geometry and hence we will have $\ttInd(\ttBan_k) = \ttInd(\ttBan_k^{nA})$. When $k$ is Archimedean, we use $\ttBan_k = \ttBan_k^{A}$. This is mainly due to state the right universal properties characterizing dagger affinoid algebras in relations to classical theories. In the Section \ref{sketch} we will see that, when working with overconvergent analytic functions, and hence with Ind-objects, there is no harm in using $\ttInd(\ttBan_k^A)$ also for non-Archimedean base fields.

\begin{defn}
Let $X$ be a set. A \emph{bornology} on $X$ is a collection $\mB$ of subsets of $X$ such that
\begin{enumerate}
\item $\mB$ is a covering of $X$, \ie $\forall x \in X, \exists B \in \mB$ such that $x \in \mB$;
\item $\mB$ is stable under inclusions, \ie $A \subset B \in \mB \then A \in \mB$;
\item $\mB$ is stable under finite unions, \ie for each $n \in \Nb$ and $B_1, ..., B_n \in \mB$, $\bigcup_{i = 1}^n B_i \in \mB$.
\end{enumerate}

The pair $(X, \mB)$ is called a \emph{bornological set}, and the elements of $\mB$ are called \emph{bounded subsets} of $X$ 
(with respect to $\mB$, if it is needed to specify).  
A family of subsets $\A \subset \mB$ is called a \emph{basis} for $\mB$ if for any $B \in \mB$ there exist $A_1, \dots, A_n \in \A$ such that $B \subset A_1 \cup \dots \cup A_n$.
A \emph{morphism} of bornological sets $\varphi: (X, \mB_X) \to (Y, \mB_Y)$ is defined to be a bounded map $\varphi: X \to Y$, \ie a map of sets such that 
$\varphi(B) \in \mB_Y$ for all $B \in \mB_X$. 
\end{defn}

The category of bornological sets is complete and cocomplete. For a detailed study of general facts about bornological sets and bornological algebraic structures we refer to \cite{Bam},
chapter 1. The category of bornological sets has some analogies with the category of topological spaces, for example the forgetful functor from bornological sets to sets is a
topological functor (in the sense of \cite{ACC}), but the theory of bornological sets seems to don't have a direct geometrical interest as the category of topological spaces.
From now on we focus on the notion of bornological vector spaces over a complete non-trivially valued field $k$ (Archimedean or non-Archimedean). The valuation on $k$ defines
a notion of boundedness with respect to which $k$ is a bornological ring.

\begin{defn}
A \emph{bornological vector space} over $k$ is a $k$-vector space $E$ along with a bornology on the underlying set of $E$ for which the maps $(\la, x) \mapsto \la x$ and $(x, y) \mapsto x + y$
are bounded.
\end{defn}

We will be particularly interested on bornological vector spaces whose bounded subsets can be described using convex subsets, in the following way.

\begin{defn}
Let $E$ be a $k$-vector space. A subset $B \subset E$ is called \emph{absolutely convex} (or a \emph{disk}) if
\begin{enumerate}
\item for $k$ Archimedean, it is convex and balanced, where convex means that for every $x, y \in B$ and $t \in [0, 1]$ then $(1 -t) x + t y \in B$ and balanced means that for any $r \in k^\circ$, $r B \subset B$;
\item for $k$ non-Archimedean, it is a $k^\circ$-submodule of $E$.
\end{enumerate}
\end{defn}

The definition of absolutely convex subset of $E$ is posed in two different ways, depending on $k$ being Archimedean of non-Archimedean, although the formal properties are the same in both cases. One can give a uniform description of absolutely convex subsets of $k$-vector spaces using generalized rings, for example as done by Durov in \cite{Dur}.
 
\begin{defn}
A bornological vector space is said to be \emph{of convex type} if it has a basis made of absolutely convex subsets. We will denote by $\ttBorn_k$ the category whose objects are the  bornological vector spaces of convex type and whose morphisms are bounded linear maps between them.
\end{defn}

\begin{rem}
For every bornological vector space of convex type $E$ there is an isomorphism
\[ E \cong \limind_{B \in \mB_E} E_B  \]
where $B$ varies over the family of bounded absolutely convex subsets of $E$ and $E_B$ is the vector subspace of $E$ spanned by elements of $B$  equipped with the gauge semi-norm (also called Minkowski functional) defined by $B$.
\end{rem}

\begin{defn} \label{defn:separated_born}
        A bornological vector space over $k$ is said to be \emph{separated} if its only bounded vector subspace is the trivial subspace $\{0\}$.
\end{defn}
\begin{rem}
A bornological vector space of convex type over $k$ is separated if and only if for each $B\in \mathcal{B}_{E}$, the gauge semi-norm on $E_{B}$ is actually a norm.
\end{rem}
The category of separated bornological spaces of convex type over $k$ is denoted $\ttSBorn_k$ and is fully faithfully embedded in $\ttBorn_k$.

\begin{defn} \label{defn:complete_born}
        A bornological space $E$ over $k$ is said to be \emph{complete} if there is a small filtered category $I$, a functor $I \to \ttBan_{k}$ and an isomorphism
        \[ E \cong \limind_{i \in I} E_i \]
        for a filtered colimit of Banach spaces over $k$ for which the system morphisms are all injective. Note that this is a colimit in the category $\ttBorn_k$.
\end{defn}

The category of complete bornological spaces over $k$ is denoted $\ttCBorn_k$ and is by definition  fully faithfully embedded in $\ttSBorn_k$, obtaining the following fully faithful embeddings:
\[ \ttCBorn_k \rhook \ttSBorn_k \rhook \ttBorn_k. \]

The inclusion functor $\ttSBorn_k \rhook \ttBorn_k$ admits a left adjoint which is given in the following way: let $E \in \ttBorn_k$, the separation of $E$ is given by
\[ \sep (E) = \frac{E}{\ol{\{0\}}}. \]
where $\ol{\{0\}}$ means the bornological closure of $\{0\}$ (see below for the definition of bornological closure). Also the inclusion functor $\ttCBorn_k \rhook \ttSBorn_k$ admits a left adjoint which is given in the following way: let $E \in \ttSBorn_k$, the completion of $E$ is given by
\[ \what{E} = \limind_{B\in \mathcal{B}_{E}} \what{E}_B. \]
where $\limind_{B\in \mathcal{B}_{E}}$ is calculated in $\ttSBorn_{k}$. 

The fact that this completion functor $\ttSBorn_{k} \to \ttCBorn_{k}$ is left adjoint to the inclusion follows from the computation that for any $E\in \ttSBorn_{k}$ and $F \in \ttCBorn_{k}$ 
\[
\begin{split} \Hom_{\ttCBorn_{k}}(\widehat{E},  F ) & = \Hom_{\ttCBorn_{k}}(\limind_{B \in \mathcal{B}_{E}} \widehat{E}_B,  \limind_{i \in I} F_i ) =  \limpro_{B \in \mathcal{B}_{E}} \limind_{i\in I}  \Hom_{\ttBan_{k}}(\widehat{ E}_B, F_i)   \\ & =   \limpro_{B\in \mathcal{B}_{E}} \limind_{i\in I}  \Hom_{\ttNrm_{k}}( E_B, F_i) =\Hom_{\ttSBorn_{k}}(\limind_{B \in \mathcal{B}_{E}} E_B,  \limind_{i \in I} F_i ) \\ & = \Hom_{\ttSBorn_{k}}(E,  F ) 
\end{split}
\]
The left adjointness of the separation functor is entirely analogous. Composing the two adjoint functors to the inclusions, one can define the (separated) completion of any bornological vector space of convex type. For $E\in \ttBorn_{k}$ we often write $\widehat{E}$ in place of $\widehat{\sep(E)}$.

\begin{rem}
Limits and colimits in $\ttBorn_{k}$ and $\ttCBorn_{k}$ are easy to describe explicitly.
Suppose we are given a functor $F:S \to \ttBorn_{k}$ where $S$ is a small category. Let $V_{F}$ be the limit of $F$ in the category of $k$-modules. Declare a subset of $V_{F}$ to be bounded if its image in each $F(s)$ is bounded. This defines a bornology on $V_{F}$.  The morphism $V_{F} \to F(s)$ is then clearly bounded for each $s$. Given a bornological vector space $V$ over $k$ and a compatible set of bounded linear maps $V \to F(s)$ consider the linear map $V \to V_{F}$ coming from the universal property of limits in the category of vector spaces. The image of a bounded set in $V$ inside $V_{F}$ is bounded by the factorization structure. For the colimits, let $V^{F}$ be the colimit of $F$ in the category of $k$-modules. Consider the vector space bornology on $V^F$ generated by the images of bounded sets from the $F(s)$.  Given a bornological vector space $V$ over $k$ and a compatible set of bounded linear maps $F(s) \to V $ consider the map $V^{F} \to V$ coming from the universal property of colimits in the category of vector spaces. This map sends the generating sets to bounded sets and so it is bounded. For a functor, $F:S \to \ttCBorn_{k}$, the limit  of the postcomposition into $\ttBorn_{k}$ is complete and gives the limit in  $\ttCBorn_{k}$. For colimits, one needs only to do the separated completion of the postcomposition.
\end{rem}

The equiboundedness bornology on $\Hom_{\ttBorn_k}(E, F)$ is a vector space bornology which induces an internal $\Hom$ functor denoted $\uHom_{\ttBorn_k}$. 
Moreover, given two bornological vector spaces $E, F$ of convex type one can put on $E \otimes_{k} F$ the \emph{projective tensor product bornology} in the following way:
a basis for the bornology of $E \otimes_{\pi,k} F$ is given by absolutely convex hulls of subsets of the form
\[ A \otimes B = \{ x \otimes y \ \ | \ \ x \in A, y \in B \} \]
where $A, B$ varies over a basis of absolutely convex subsets for the bornologies of $E$ and $F$.
We have functors 
\begin{equation}\label{eqn:Bigger}
\begin{split}
& \ttSNrm_{k} \to \ttBorn_k \\
& \ttNrm_{k} \to \ttSBorn_k \\
& \ttBan_{k} \to \ttCBorn_k
\end{split}
\end{equation}
where the bornology on the underlying vector space of a semi-normed space is given by the bounded subsets with respect to the semi-norm. The following is a trivial but important observation.

\begin{lem} \label{lemma:born_elementary}
The category $\ttCBorn_k$ is elementary quasi-abelian with all limits and colimits and enough flat projectives.
\end{lem}

{\bf Proof.}
A proof of this can be found in \cite{PrSc}, where only the case $k = \Cb$ is treated. But in fact the same argument of \cite{PrSc} can be used for any valued field once one define coproducts in the contracting categories $\ttBan_k^{\le 1}$ as we did in Proposition \ref{prop:contracting_limits}. The construction of enough projectives works exactly as in the case of $\ttInd(\ttBan^{A}_{k})$ discussed in Lemma \ref{lem:exact_direct_limit}.

\ \hfill $\Box$

\begin{prop}
Let $E, F, G \in \ttBorn_k$ then
\[ \Hom_{\ttBorn_k}(E \otimes_{\pi,k} F, G) \cong \Hom_{\ttBorn_k}(E, \uHom_{\ttBorn_k}(F, G)). \]
\end{prop}

{\bf Proof.}
This is a restatement of theorem 1b of page 174 of \cite{H2}, where in fact is proved that 
\[  \uHom_{\ttBorn_k}(E \otimes_{\pi,k} F, G) \cong \uHom_{\ttBorn_k}(E, \uHom_{\ttBorn_k}(F, G)). \]
\ \hfill $\Box$ 

From the presentation $E \cong \limind_{B \in \mB_E} E_B$ for any object of $\ttBorn_k$ we can define the association 
\[E \mapsto ``\limind_{B \in \mB_E}" E_B
\] which extends in a clear way to functors, sometimes called \emph{dissection functors}, 
\[ \diss: \ttBorn_k \to \ttInd(\ttSNrm_k) \]
\[ \diss: \ttSBorn_k \to \ttInd(\ttNrm_k) \]
\[ \diss: \ttCBorn_k \to \ttInd(\ttBan_k). \]
They are fully faithful and their essential image is given by the essentially monomorphic objects, see Definition \ref{defn:essential_monomorphic}. 

As done in previous sections we will focus on the study of the category of complete objects, hence on studying $\ttCBorn_k$ and its relations with $\ttInd(\ttBan_k)$.
We now describe the monoidal structure on $\ttCBorn_k$ which is known as the \emph{completed projective tensor product}. It is defined by  
\begin{defn}\label{defn:BornModStr}
\[ E \wotimes_{\pi,k} F = \what{E \otimes_{\pi,k} F} \]
for any $E, F \in \ttCBorn_k$. Here we are implicitly using the fact that $E \otimes_{\pi,k} F$ is always separated for $E$ and $F$ separated. Unless there is some other monoidal structure being discussed, we simplify the notation, writing $E \wotimes_{k} F$ in place of $E \wotimes_{\pi,k} F$.
\end{defn}
\begin{rem}\label{rem:CommBorn}
With the monodial structure from Definition \ref{defn:BornModStr} we see that $\ttCBorn_{k}$ is a closed, symmetric monoidal category:
\[  \uHom_{\ttCBorn_k}(E \wotimes_{\pi,k} F, G) \cong \uHom_{\ttCBorn_k}(E, \uHom_{\ttCBorn_k}(F, G)). \]
Therefore, we can consider the category $\ttComm(\ttCBorn_k)$. This category is convenient to work with because in \cite{Bam} it is shown that the forgetful functor from $\ttComm(\ttCBorn_k) \to \ttCBorn_k$ commutes with all limits and also with filtered colimits.
\end{rem}
\begin{obs}The functors (\ref{eqn:Bigger}) are fully faithful, respect strict morphisms, monomorphisms and epimorphisms, finite limits and finite colimits the internal Hom and monoidal structures. The image of any projective (respectively flat) object under any of these functors is projective (respectively flat). In particular, a morphism $A \to B$ in $\ttComm(\ttBan_{k})$ is a homotopy epimorphism if and only if it becomes one in $\ttComm(\ttCBorn_{k})$
\end{obs}
\begin{prop} 
The dissection functor $\diss: \ttCBorn_k \to \ttInd(\ttBan_k)$:
\begin{enumerate}
\item define an equivalence of $\ttCBorn_k$ with the sub-category of essential monomorphic objects of $\ttInd(\ttBan_k)$;
\item commutes with all limits and coproducts;
\item preserves the internal hom;
\item does not commute with filtered colimits or cokernels in general;
\item commutes with monomorphic filtered colimits. 
\end{enumerate} \end{prop}

{\bf Proof.} For the first three properties see \cite{M}, Theorem 1.139.  To see that $\diss$ doesn't commute with filtered colimits is enough to consider and element which is not in the essential image of $\diss$ and calculate its colimit in $\ttCBorn_k$. For a counter-example with cokernels one can do an analogous reasoning as done in \ref{rem:non_concrete_ban}. Finally, $\diss$ commutes with monomorphic filtered colimits because the composition of two filtered monomorphic colimits is a monomorphic filtered colimit.

\ \hfill $\Box$

\begin{rem}
The dissection functor $\diss: \ttCBorn_k \to \ttInd(\ttBan_k)$ unfortunately does not respect the monoidal structures. The problem with this functor lies on the left hand side, in the notion of complete bornological vector space. Namely, this notion is not so good as in the case of topological vector spaces and presents some anomalies: for example the completion of a separated bornological vector space can be the null space. To avoid this ill-behaviour one can introduce the notion of properness. To do this we need to introduce the notion of bornological (or Mackey) convergence and bornological closed subsets.
\end{rem}

\begin{defn}
Let $E$ be a bornological $k$-vector space and $\{x_n\}$ a sequence of elements of $E$. We say that $\{x_n\}$ \emph{converges (bornologically) to $0$ in the sense of Mackey} if there exists a bounded subset $B \subset E$ such that for every $\la \in k^\times$ there exists an $n = n(\la)$ for which
\[ \{x_m\} \subset \la B, \forall m > n. \] 
We say that $\{ x_n \}$ converges  (bornologically) to $a \in E$ if $\{x_n - a\}$ converges (bornologically) to zero.
\end{defn}

In analogous way one can give the definition of convergence of filters of subsets of $E$. We omit the details of this definition since is not important for our scope.

\begin{defn}
Let $E$ be a bornological vector space over $k$.
\begin{itemize}
\item a sequence $\{x_n\} \subset E$ is called \emph{Cauchy-Mackey} if the double sequence $\{x_n - x_m\}$ converges to zero;
\item a subset $U \subset E$ is called \emph{(bornologically) closed} if every Mackey convergent sequence of elements of $U$ converges (bornologically) to an element of $U$. 
\end{itemize}
\end{defn}

\begin{defn}
A bornological vector space is called \emph{semi-complete} if every Cauchy-Mackey sequence is convergent.
\end{defn}

The notion of semi-completeness is not as useful as the notion of completeness in the theory of topological vector spaces. We remark that any complete bornological vector space is semi-complete, but the converse is false.

\begin{rem}
The notion of bornological convergence of bornological vector spaces of convex type $E = \limind_{B \in \mB_E} E_B$, where $B$ varies over the family of bounded disks of $E$, can be restated in the following way: $\{ x_n \}_{n \in \Nb}$ is convergent to zero in the sense of Mackey if and only if there exists an $B \in \mB_E$ and and $N \in \Nb$ such that for all $n > N$, $x_n \in E_B$ and $x_n \to 0$ in $E_B$ for the semi-norm of $E_B$.
\end{rem}

It can be shown that the notion of bornologically closed subset induces a topology on $E$, but this topology in neither a vector space topology nor group topology in general. Therefore, an arbitrary intersection of bornological closed subsets of a bornological vector space is bornologically closed. So, the following definition is well posed.

\begin{defn}
 Let $U \subset E$ be a subset of a bornolgical vector space. The closure of $U$ is the smallest bornologically closed subset of $E$ in which $U$ is contained. We denote the closure of $U$ by $\ol{U}$.
\end{defn}

The concept of bornologically closed subspace fits nicely in the theory. For example a bornological vector space is separated, in the sense of Definition \ref{defn:separated_born}, if and only if $\{0 \}$ is a bornologically closed subset.

\begin{defn} \label{defn:born_dense}
Let $E$ be a bornological vector space over $k$. We will say that a subset $U \subset E$ is \emph{bornologically dense} if the bornological closure of $U$ is equal to $E$. 
\end{defn}
\begin{obs}\label{obs:Dense2Dense}Let $E$ be a bornological vector space over $k$. If $U \subset E$ is bornologically dense then $U \subset E^{t}$ is dense. 
\end{obs}

The following proposition has an analog for $\ttBorn_k$ and $\ttSBorn_k$, but we state it only for $\ttCBorn_k$ since is the only case we are interested in.

\begin{prop} \label{prop:strict_morphisms_born_close}
Let $f: E \to F$ be a morphism in $\ttCBorn_k$, then
\begin{itemize}
\item $f$ is a monomorphism if it is injective
\item $f$ is an epimorphism if $f(E)$ is bornologically dense in $F$
\item $f$ is a strict epimorphism if and only if it is surjective and $F$ is endowed with the quotient bornology; 
\item $f$ is a strict monomorphism if and only if it is injective, the bornology on $E$ agrees with the induced bornology from $F$ and is a closed subspace of $F$.
\end{itemize}
\end{prop}

{\bf Proof.}

It is easy to deduce these claims from what we discussed so far, and we omit the details that can be found in \cite{PrSc} or \cite{H2}.

\ \hfill $\Box$

After this recall we can state the notion of properness.

\begin{defn} \label{def:proper_born}
A bornological vector space is called \emph{proper} if its bornology has a basis of bornologically closed subsets.
\end{defn}

\begin{rem} \label{rem:proper}
All spaces considered in \cite{Bam} satisfies this property. One can show that a separated proper bornological vector space injects in its completion and moreover the following Proposition.
\end{rem}

\begin{prop} \label{prop:proper}
The dissection functor from the full sub-category of proper objects in $\ttCBorn_k$ to $\ttInd(\ttBan_k)$ respects the monoidal structures.
\end{prop}

{\bf Proof.}

Proposition 1.149 and Corollary 1.151, of \cite{M} gives the assertion. The key point is that for a proper bornological vector spaces $E$ the map $E \to \widehat{E}$ is injective, \ie the space $\limind_{B\in \mathcal{B}_{E}} \what{E}_B$ is a separated, complete (and also proper) bornological vector space.

\ \hfill $\Box$

\section{Dagger affinoid algebras and spaces}

Here we recall definitions and properties of the theory of dagger affinoid algebras as developed in \cite{Bam}, which will be the main reference for this section. Moreover we will, without any harm to final results, recast all the statements in the category of $\ttInd(\ttBan_k)$ instead of $\ttCBorn_k$, in the view of pursuing the point of view that will be discussed with more details Section \ref{sketch}. In this section, $k$ will be still supposed to be a complete non-trivially valued field.

\begin{defn}
We say that an algebra $A \in \ttComm(\ttInd(\ttBan_k))$ is \emph{multiplicatively convex} (or an \emph{m-algebra}) if $A$ is isomorphic to an object of $\ttInd(\ttComm(\ttBan_k))$.
\end{defn}

Hence, the difference between an algebra on $\ttInd(\ttBan_k)$ and a multiplicatively convex algebra over $k$ is that the former is an algebra object that can be presented as a directed system of Banach spaces, while the latter can be presented as the colimit of a directed system Banach algebras.
These two categories do not agree, for explicit counter-examples see \cite{Bam} and \cite{M}.

\begin{defn} \label{defn:Berkovich_spectrum}
Let $A \in \ttInd(\ttComm(\ttBan_k))$ we define the \emph{(Berkovich) spectrum}
of $A$, denoted $\mM(A)$ as the topological space
\[ \mM(A) = \mM(``\limind_{i \in I}" A_i) \cong \limpro_{i \in I} \mM(A_i), \]
where $\mM(A_i)$ are the usual Berkovich spectra of $A_i$. $\mM(A)$ is equipped with the projective limit topology.
\end{defn}

\begin{lem} 
For $A \in \ttInd(\ttComm(\ttBan_k))$, $\mM(A)$ is a non-empty compact, Hausdorff topological space.
\end{lem}
{\bf Proof.}
The same argument of \cite{Bam} used for the bornological spectrum applies.
\ \hfill $\Box$

\begin{defn}
For any polyradius $\rho = (\rho_1, ..., \rho_n) \in \Rb_+^n$ we define an object of $\ttComm(\ttInd(\ttBan_k))$ by 
\[ \mW_k^n (\rho) = k \lt \rho_1^{-1} X_1, ..., \rho_n^{-1} X_n \gt^\dagger = ``\limind_{r > \rho}" \mT_k^n(r) \]
where $\mT_k^n(r)$ is the ring of ``strictly convergent" power-series in the polycylinder 
\[ \{c \in k^{n} \ \ | \ \ |c_i| \leq r_{i} \}\] 
of polyradius $r$, and call it the ring of \emph{overconvergent analytic functions} on the polycylinder of polyradius $\rho = (\rho_1, ..., \rho_n)$ centred in zero. More explicitly, if $k$ is non-Archimedean then 
\[\mT_k^n(r) = \{ \sum_{I \in \Nb^n} a_I X^I \ \  | \ \ \lim_{I\to \infty} |a_I| r^I  \to 0 \}, \]
equipped with the norm $\|\sum_{I \in \Nb^n} a_I X^I\| = \sup_{I \in \Nb^n}|a_I| r^I$
is a non-Archimedean Banach algebra called the Tate algebra of the polycylinder of polyradius $r$. If $k$ is Archimedean then
\[ \mT_k^n(r) = \{ \sum_{I \in \Nb^n} a_I X^I \ \ | \ \ \sum_{I \in \Nb^n}|a_I| r^I < \infty \}, \]
equipped with the norm $\|\sum_{I \in \Nb^n} a_I X^I\| = \sum_{I \in \Nb^n}|a_I| r^I$ is an Archimedean Banach algebra.
Moreover, if $\rho = (1, ..., 1)$ we simply write
\[ \mW_k^n = k \lt X_1, ..., X_n \gt^\dagger \]
and call it the ring of \emph{overconvergent analytic functions} on the polydisk of radius $1$.
\end{defn}

In \cite{Bam} we consider on $\mW_k^n (\rho)$ the direct limit bornology induced by the $k$-Banach algebra structures on $\mT_k^n(\rho)$, here we consider $\mW_k^n (\rho)$ as an algebra object on $\ttInd(\ttBan_k)$. By the discussion of previous section is clear that this two point of view are essentially the same.

\begin{rem}
$\mW_k^n(\rho)$ could also be defined as the direct limit of the Frech\'{e}t algebras of analytic functions on open polycylinders of radius bigger than $\rho$. If $k$ is Archimedean, it could also be defined as the direct limit of disk algebras (we will recall the definition of disk algebras in the Section \ref{sketch} of the article). It is a non-trivial result that all these definitions of $\mW_k^n(\rho)$ are equivalent. We will say more on this issue in Section \ref{sketch} of the paper.
\end{rem}

We defined $\mW_k^n(\rho)$ as an object of $\ttComm(\ttInd(\ttBan_k))$ and the underlying $\ttInd(\ttBan_k)$ object of $\mW_k^n(\rho)$ is essentially monomorphic. So by Corollary \ref{prop:essential_concrete} we know that is meaningful to associate to $\mW_k^n(\rho)$ its underlying ring given by
\[ \bigcup_{r > \rho} \mT_k^n(r), \]
Moreover, by Proposition \ref{prop:essentially_subobject} we also know that all ideals (sub-objects in the category of non-unital algebras relative to $\ttInd(\ttBan_k)$) of $\mW_k^n(\rho)$, as object of $\ttComm(\ttInd(\ttBan_k))$, corresponds to ideals of its underlying ring. So, there will be no confusion on what we mean by an ideal in the following.

\begin{rem}
We don't discuss here the case when $k$ is trivially valued because in this case there is no relation between bornological spaces and $\ttInd(\ttBan_k)$. This case is more similar to the case of a general base Banach ring and an approach in this direction is sketched in Section \ref{sketch}, where a formal approach is proposed to overcome difficulties. 
\end{rem}

\begin{defn}\label{defn:AffinoidAlg}
An algebra $\A \in \ttComm(\ttInd(\ttBan_k))$ is called a \emph{dagger affinoid algebra} if 
\[ \A \cong \frac{\mW_k^n(\rho)}{I} \]
for some $n \in \Nb$, $\rho \in \Rb_+^n$ and ideal $I \subset \mW_k^n(\rho)$. We denote by 
\[ \ttAfnd_k^\dagger \rhook \ttComm(\ttInd(\ttBan_k)) \] the full sub-category identified by dagger affinoid algebras.
\end{defn}

We sum up in the next theorem the main properties of dagger affinoid algebras. We remark that by definition the underlying $\ttInd(\ttBan_k)$ object of a dagger affinoid algebra is essentially monomorphic, hence we can associate to any dagger affinoid algebra a meaningful underlying ring, and we tacitly refer to this ring when we will talk about ring-theoretic properties of dagger affinoid algebras.

\begin{thm}  \label{thm:dagger_algebras}
\begin{enumerate}
\item All objects of $\ttAfnd_k^\dagger$ are Noetherian and $\mW_k^n (\rho)$ is factorial.
\item All the ideals of objects of $\ttAfnd_k^\dagger$ are closed, \ie the operation of quotienting by ideals, calculated in $\ttInd(\ttSNrm_k)$, gives separated objects \ie objects in the sub-category $\ttInd(\ttBan_k)$.
\item The forgetful functor $\ttAfnd_k^\dagger \to \ttComm(\ttVect_k)$ is fully faithful;
\item Every morphism in $\ttAfnd_k^\dagger$ is a morphism of inductive systems of $k$-Banach algebras, \ie the embedding of $\ttAfnd_k^\dagger \to \ttInd(\ttComm(\ttBan_k))$ is fully faithful.
\end{enumerate}\end{thm}

{\bf Proof.}

The proofs of this results can be found in the first two sections of chapter 3 of \cite{Bam}.

\ \hfill $\Box$

\begin{defn}\label{defn:DagAffLoc}
A \emph{$k$-dagger affinoid localization} is a morphism $\A \to \D$ of $k$-dagger affinoid algebras such that the morphism $\mM(\D) \to \mM(\A)$ is injective and any morphism of $k$-dagger affinoid algebras $\A \to \B$ such that $\mM(\B)$ lands in the image of $\mM(\D)$ factors as $\A \to \D \to \B$. The subset $V \subset \mM(\A)$ identified by a $k$-dagger affinoid localization is called a \emph{dagger affinoid subdomain} of $\mM(\A)$ and we often write $D=A_{V}$. The morphism $V=\mM(D) \to \mM(\A)$ is called a \emph{$k$-dagger affinoid immersion}.
\end{defn}

The family of dagger affinoid subdomains satisfies the property of a pre-topology, and the following particular $k$-dagger affinoid localizations are the ones used in calculations.

\begin{defn}
A \emph{Weierstrass localization} of a $k$-dagger affinoid algebra $\A$ is a morphism of the form
\[ \A \to \frac{ \A \lt r_1^{-1} X_1, \ldots, r_n^{-1} X_n \gt^\dagger}{(X_1 - f_1, \ldots, X_n - f_n)} \]
for some $f_1, \dots, f_n \in \A$, $(r_i) \in \Rb_+^n$. A \emph{Laurent localization} of a $k$-dagger affinoid algebra $\A$ is a morphism of the form
\[ \A \to \frac{ \A \lt r_1^{-1} X_1, \ldots, r_n^{-1} X_n, s_1^{-1} Y_1, \ldots, s_m^{-1} Y_m \gt^\dagger}{(X_1 - f_1, \ldots, X_n - f_n, g_1 Y_1 - 1, \ldots, g_m Y_m - 1)} \] 
for some $f_1, \dots, f_n, g_1, \dots, g_m \in \A$, $(r_i), (s_i) \in \Rb_+^n$. Finally, a \emph{rational localization} of a $k$-dagger affinoid algebra $A$ is a morphism of the form
\[ \A \to \frac{ \A \lt r_1^{-1} X_1, \ldots, r_n^{-1} X_n \gt^\dagger}{(h X_1 - f_1, \ldots, h X_n - f_n)} \] 
for some $f_1, \dots, f_n, h \in \A$, $(r_i) \in \Rb_+^n$ such that the identity of $A$ belongs to the ideal generated by $h, f_{1}, \dots, f_{n}$. The corresponding maps $\mathcal{M}(A_V) \to \mathcal{M}(A)$ are called Weierstrass immersions, Laurent immersions and rational immersions.
\end{defn}

We assume that the reader is familiar with the category of locally ringed Grothendieck topological spaces, as explained in \cite{BGR}, chapter 9. 
\begin{defn} \label{defn:stalks}
Let $X$ be a locally ringed Grothendieck topological space and $x \in X$. The stalk of $\mO_X$ at $x$ is defined to be
\[ \mO_{X, x} = \limind \mO_X(U) \]
where $U$ runs over the family of admissible open subsets that contains $x$. Notice that we are treating all the objects here as belonging to the category of rings (with no extra structure).
\end{defn}
\begin{defn}\label{defn:openImmersionGRS}
An open immersion $X\to Y$ between locally ringed Grothendieck topological spaces is a morphism in the category of locally ringed Grothenidieck topological spaces which is injective on sets and an isomorphism on all stalks.
\end{defn}

\begin{defn}
The Grothendieck topological space whose underlying set is the spectrum of a $k$-dagger affinoid algebra $A$ and whose covering families are the finite covers by $k$-dagger affinoid subdomains is called a \emph{$k$-dagger affinoid space}.
\end{defn}

The Grothendieck described in previous definition is called the \emph{weak dagger G-topology} and $\mM(A)$ equipped with this topology is denoted $\mM^G(A)$.

\begin{defn}\label{defn:DaggerAffinoidSpaces} The category of dagger affinoid spaces over $k$ is defined to be the full sub-category of locally ringed Grothendieck topological spaces of the form $\mM^{G}(A)$ where $A$ is a dagger affinoid $k$-algebras where the structure sheaf is defined by $V \mapsto A_{V}$ (explained in Definition \ref{defn:DagAffLoc}).
\end{defn}

\begin{rem}
One can check using Theorem \ref{thm:dagger_algebras} that the category of dagger affinoid algebras and the category of dagger affinoid spaces are anti-equivalent.
\end{rem}
\begin{rem}\label{rem:usual}
        The usual properties of classical affinoid immersions hold for dagger affinoid immersions. We list here the main ones, and we refer to the last part of chapter 3 of \cite{Bam} for proofs:
        \begin{enumerate}
                \item if $U \rhook V$ is a dagger affinoid immersion and $V \rhook X$ is another one, then the composition $U \rhook X$ is a dagger affinoid immersion;
                \item if $U \subset X$ is a dagger affinoid immersion and $\phi: Y \to X$ a morphism of dagger affinoid spaces then the natural morphism $\phi^{-1}(U)\to Y$ is a dagger affinoid immersion of $Y$ and the type of the immersion (Weierstrass, Laurent, rational) is preserved;
                \item The intersection of two dagger affinoid immersions is a dagger affinoid immersion. If both immersions are Weierstrass (respectively Laurent or rational) then there intersection is Weierstrass (respectively Laurent or rational) \end{enumerate}
\end{rem}

\ \hfill $\Box$

\begin{defn}
A morphism of $k$-dagger affinoid spaces $f: X \to Y$, is called \emph{Runge immersion} if it factors in a diagram
\[
\begin{tikzpicture}
\matrix(m)[matrix of math nodes,
row sep=2.6em, column sep=2.8em,
text height=1.5ex, text depth=0.25ex]
{X & & Y  \\
& Y' \\};
\path[->,font=\scriptsize]
(m-1-1) edge node[auto] {$f$} (m-1-3);
\path[->,font=\scriptsize]
(m-1-1) edge node[auto] {$g$} (m-2-2);
\path[->,font=\scriptsize]
(m-2-2) edge node[auto] {$h$}  (m-1-3);
\end{tikzpicture}
\]
where $g$ is a closed immersion and $h: Y' \to Y$ is a Weierstrass domain immersion.
\end{defn}

The following are the main results in the theory of dagger affinoid spaces.

\begin{thm}\label{thm:MainResults}
\begin{itemize}
\item (Gerritzen-Grauert theorem) Let $\phi: X \to Y$ be a locally closed immersion of $k$-dagger affinoid spaces, then there exists a finite covering of $\phi(X)$ by rational subdomains $Y_i$ of $Y$ such that all morphisms $\phi_i: \phi^{-1}(Y_i) \to Y_i$ are Runge immersions.
\item Main corollary to Gerritzen-Grauert theorem: Let $X$ be a $k$-dagger affinoid space and $U \subset X$ a $k$-dagger affinoid subdomain then there exist a finite number of rational subdomains $U_i \subset X$ such that $\bigcup U_i = U$.
\item (Tate's acyclicity theorem) Let $A$ be a $k$-dagger affinoid algebra, then the presheaf $U \mapsto \A_U$ is acyclic for $\mM^G(\A)$.
\item (Kiehl's theorem) Let $\A$ be a $k$-dagger affinoid algebra and $X = \mM^G(\A)$. Then an $\Oc_X$-modules is coherent if and only if is associated to a finite $\A$-module.
\end{itemize}\end{thm}

{\bf Proof.} 

The proofs can be found in \cite{Bam}, chapter 4, 5 and 6.

\ \hfill $\Box$

Notice that by Tate acyclicity the presheaf $U \mapsto \A_U$, for dagger affinoid localizations is a sheaf. 

From now on if $X = \mM(\A)$ is a dagger affinoid space and $U \subset X$ is a dagger 
affinoid subdomain we denote $\A_U$ by $\mO_X(U)$ and we endow $\mM^G(\A)$ with the 
locally $G$-ringed space structure $(\mM^G(\A), \mO_X)$. 

\begin{rem}
In the case $k = \Cb$ any dagger affinoid space $X$ has a canonical structure of compact Stein subset of $\Cb^n$, for some $n$ and  all constructions above are compatible with this structure. In particular in this case the weak $G$-topology defines a sub-site of the compact Stein site of $X$. 
\end{rem}
To globalize the construction we gave so far we need to use Berkovich nets, since our building blocks are compact spaces. This procedure is long to describe in details and is totally analogous to what Berkovich did in \cite{Ber1993}. We refer to the last chapter of \cite{Bam} for details of the dagger case of this constructions. 
\begin{defn}
A \emph{$k$-dagger analytic space} is the data of a triple $(X, A, \tau)$ where $X$ is a topological space, $\tau$ is a Berkovich net on $X$ and $A$ is an atlas of dagger affinoid subdomains for $\tau$.
\end{defn}
\begin{rem}\label{rem:HasSheaf}
Any $k$-dagger analytic space $X$ carries a sheaf of bornological algebras $\mathcal{O}_{X}$ and in the case that $X=\mathcal{M}^{G}(A)$ is a $k$-dagger affinoid then $\mathcal{O}_{X}(X) \cong A$ as bornological algebras, as explained in chapter 6 of \cite{Bam}. In general, $X \mapsto \mathcal{O}_{X}(X)$ is a functor from $k$-dagger analytic spaces to bornological algebras which is right adjoint to $\mathcal{M}^{G}$.
\end{rem}
We end this section with some definitions involving stalks. Notice that a dagger affinoid subdomain $U \subset X$ might be not a neighborhood of $x \in U$ in the topology of $\mM(A)$, but this is not a problem since the elements of $\mO_X(U)$ are germs of analytic functions on $U$, hence each of them is defined on a neighborhood of $U$. Therefore is meaningful to apply to $\mM^G(A)$ the definition of stalk given in Definition \ref{defn:stalks} (be careful that this is not true in Berkovich geometry, where also non-overconvergent analytic functions are considered).
The next proposition shows that dagger affinoid immersions are precisely the open immersions in the category of dagger affinoid spaces.

\begin{prop}  \label{prop:openness_affinoid_immersion}
A morphism of $k$-dagger affinoids is a $k$-dagger affinoid immersion (see Definition \ref{defn:DagAffLoc}) if and only if it induces an open immersion in the sense of Definition \ref{defn:openImmersionGRS} of the associated locally ringed Grothendieck topological spaces.\end{prop}

{\bf Proof.}

If $f$ is dagger affinoid immersion then the isomorphism of stalks is immediate from Definition \ref{defn:stalks} and the injectivity of $f$ is one of the basic properties of dagger affinoid immersions (or imposed by definition as in this article, depending on different equivalent characterizations). On the other hand, if $f$ is an open immersion then $f$ is an isomorphism onto its image, which by the main corollary of the Gerritzen-Grauert theorem (see Theroem \ref{thm:MainResults}) is a union of $k$-dagger rational subdomains. So, by Proposition 6.1.27 of \cite{Bam} $f$ is a $k$-dagger affinoid immersion.

\ \hfill $\Box$

\section{Dagger analytic geometry}

The aim of this section is to derive a categorical characterization of the dagger weak $G$-topology defined in the previous section, in the category of dagger affinoid spaces. This characterization is analogous to the one given in \cite{BeKr} for the weak $G$-topology of classical affinoid spaces over non-Archimedean base fields. In this section we deal with the category $\ttComm(\ttCBorn_{k})$ instead of $\ttComm(\ttInd(\ttBan_{k}))$. This is harmless because we are always be dealing with $k$-dagger affinoid algebras, whose underlying space is normal and hence proper, and for these type of spaces the monoidal structures agree by Proposition \ref{prop:proper}.

\begin{lem}\label{lem:WeierLau}
Let $A$ be an object of $\ttAfnd_{k}^\dagger$ . Let $\A_V$ be a localization of Weierstrass or Laurent type. Let $\B$ be an $\A$-algebra which is also an object of $\ttAfnd_{k}^\dagger$.
Then, the natural morphism
\[ B \wotimes_\A^\Lb \A_V \longrightarrow B \wotimes_A A_V \]
is an isomorphism in $D^{\leq 0}(\A)$. In particular, by taking $\B = \A_V$ we see that the morphism $\A \to \A_V$ is a homotopy epimorphism.\end{lem}

{\bf Proof.}
Every Weierstrass (respectively Laurent) localization of a $k$-dagger affinoid algebra can be viewed as an iteration of basic Weierstrass (respectively Laurent) localizations where only one new variable is added at each iteration. Therefore, because the composition of homotopy epimorphisms is a homotopy epimorphism, to prove this lemma, it suffices to do the Weirstrass and Laurent cases where only one extra variable is added.
It is enough to show that any object $C$ of $\ttAfnd_{k}^\dagger$ and any $f, g \in C$ the maps
\begin{equation}\label{eqn:mono1}\phi_f: C \lt r^{-1} X \gt^\dagger  \stackrel{X-f}\longrightarrow C \lt r^{-1} X \gt^\dagger \end{equation}
\begin{equation}\label{eqn:mono2}\phi_g: C \lt r^{-1} X \gt^\dagger \stackrel{g X-1}\longrightarrow C \lt r^{-1} X \gt^\dagger \end{equation}
are strict monomorphisms in $\ttMod(C)$. In fact, specializing this claim for $C = A$ we obtain strict morphisms
\[\phi_f: A \lt r^{-1} X \gt^\dagger  \stackrel{X-f}\longrightarrow A \lt r^{-1} X \gt^\dagger, \]
\[\phi_g: A \lt r^{-1} X \gt^\dagger \stackrel{g X-1}\longrightarrow A \lt r^{-1} X \gt^\dagger \]
which define a strict resolution of $A_V$. The terms in these resolutions are projective using the fact that $A \lt r^{-1} X \gt^\dagger$ is a filtered colimit of projectives in $\ttMod(A)$ along with \ref{lem:FiltProj}. Hence they are hence also $\wotimes_A$-acyclic. Taking the completed
tensor product over $A$ with $B$ we find a representative for $B \wotimes_A A_V$ which looks like
\[\phi_f: B \lt r^{-1} X \gt^\dagger  \stackrel{X-f}\longrightarrow B \lt r^{-1} X \gt^\dagger, \]
or
\[\phi_g: B \lt r^{-1} X \gt^\dagger \stackrel{g X-1}\longrightarrow B \lt r^{-1} X \gt^\dagger \]
where the $f$ and $g$ here are the images of the original ones in $B$. But is clear that these complexes are quasi-isomorphic to $B \wotimes_A A_V$, proving the lemma. Therefore, to conclude the proof, let's prove that the morphisms in equations (\ref{eqn:mono1}) and (\ref{eqn:mono2}) are strict monomorphisms. First of all we remark that we need only to show that the morphisms are injective because $\phi_f (C \lt r^{-1} X \gt^\dagger )$ and $\phi_g(C \lt r^{-1} X \gt^\dagger )$ are ideals of $C \lt r^{-1} X \gt^\dagger$, and hence (bornologically) closed by Theorem \ref{thm:dagger_algebras} and hence, by Proposition \ref{prop:strict_morphisms_born_close}, if $\phi_f$ and $\phi_g$ are injective then they are strict monomorphisms.

The injectivity of the map $\phi_g$ follows by the following easy calculations and induction. Let's suppose $a,b \in C \lt r^{-1} X \gt^\dagger$ with 
\[  a = \sum^{\infty}_{i=0} a_i (r^{-1} X)^i, b = \sum^{\infty}_{i=0} b_i (r^{-1} X)^i \] 
and $a_0, b_0 \ne 0$ (we can always reduce to this case). Then, suppose that
\[ (g X - 1) a = (g X - 1) b \]
the comparing coefficients we have that $a_0 = b_0$, and
\[ -a_1 + g a_0 = -b_1 + g b_0 \then b_1 = a_1 \]
and by induction $a_i = b_i$ for any $i$.

We are left to prove that $\phi_f$ is injective. Let's consider the linear maps $\mu_{f^n}: C \to C$ given by $\mu_{f^n}(a) = a f^n$. Since $C$ is Noetherian by Theorem \ref{thm:dagger_algebras} (2) the ascending chain of ideals
\[ \ker(\mu_f) \subset \ker(\mu_{f^2}) \subset \ldots \]
must stabilize at some $N \in \Nb$. Hence if 
\[ (X - f) a = 0 \]
for $a \in C \lt r^{-1} X \gt^\dagger$ as before, then calculating the coefficients one gets
\[ a_N f^N = a_0 \then a_N f^{N+1} = a_0 f = 0. \]
Thus, $a_N \in \ker(\mu_{f^{N + 1}}) = \ker(\mu_{f^N})$, therefore
\[ a_N f^N = 0 = a_0, \]
but by hypothesis $a_0 \ne 0$. Hence $\phi_f$ is injective and the lemma is proved.
\ \hfill $\Box$

\begin{lem} \label{lemma_subdomain}
Let $\A \in \ttAfnd_k^\dagger$ and Let $\A_V$ be an affinoid localization of $\A$ of rational type. Let $\B$ be an $\A$-algebra which is a dagger affinoid algebra over $k$.
Then the natural morphism
\[ \B \wotimes_A^\Lb \A_V \longrightarrow \B \wotimes_\A \A_V \]
is an isomorphism in $D^{\leq 0}(A)$. In particular, by taking $\B = \A_V$ we see that the morphism $\A \to \A_V$ is a homotopy epimorphism.\end{lem}

{\bf Proof.}

We can reduce the rational case to the Laurent and Weierstrass case in the following way. Consider the rational localization
\[ A_V = \frac{A \lt r_1^{-1} X_1, \ldots, r_m^{-1} X_m \gt^\dagger}{(g X_1 - f_1, \ldots,  g X_m - f_m )} \]
where the $f_i$ together with $g$ generate the unit ideal.
The following inequality holds
\[ |g|_\sup = \max_{|\cdot| \in \mM(A_V)} |g(x)| > 0 \]
for trivial reasons. Hence there exist $\epsilon > 0$ such that
\[ \mM(A_V) \subset \mM(A_W) \subset \mM(A) \]
and 
\[ A_W = \frac{A \lt \epsilon Y \gt^\dagger}{(g Y - 1)}. \]
If $\phi: A \to A_W$ is the canonical map, then $\phi(g)$ is a unit in $A_W$ and hence
\[ A_V \cong \frac{A_W \lt r_1^{-1} X_1, \ldots, r_m^{-1} X_m \gt^\dagger}{( X_1 - \frac{\phi(f_1)}{\phi(g)}, \ldots,  X_m - \frac{\phi(f_m)}{\phi(g)} )} \]
which is a composition of a Laurent localization and a Weierstrass localization. Since the composition of two homotopy empimorphisms is a homotopy epimorphism, the lemma is proven.

\ \hfill $\Box$

\begin{lem}\label{lem:CapCup}
Let $\A_{V_1}$ and $\A_{V_2}$ be two $k$-dagger rational localization of the $k$-dagger affinoid algebra $\A$ such that $\A_{V_1 \cup V_2}$ is a $k$-dagger affinoid localization. Then, for any morphism of $k$-dagger affinoid algebras $A\to B$, the natural morphism $B\wotimes^{\mathbb{L}}_{A}\A_{V_1 \cup V_2} \to B\wotimes_{A}\A_{V_1 \cup V_2}$ is an isomorphism. Therefore by taking $B=  \A_{V_1 \cup V_2}$ we see that $\A\to \A_{V_1 \cup V_2}$ is a homotopy epimorphism. \end{lem}

{\bf Proof.}
The lemma can be deduced from the strict short exact sequence
\[ 0 \to A_{V_1 \cup V_2} \to A_{V_1} \times A_{V_2} \to A_{V_1 \cap V_2} \to 0\]
where we note that $A_{V_1 \cap V_2} = A_{V_1}\wotimes_{A}A_{V_2}$. By applying the functor $V\mapsto V\wotimes^{\mathbb{L}}_{A}B$ we get the exact triangle 
\[A_{V_1 \cup V_2} \wotimes^{\mathbb{L}}_{A}B\to (A_{V_1} \wotimes^{\mathbb{L}}_{A}B)\times (A_{V_2} \wotimes^{\mathbb{L}}_{A}B) \to A_{V_1 \cap V_2}\wotimes^{\mathbb{L}}_{A}B
\]

which by Lemma \ref{lemma_subdomain} reduces to 
\[A_{V_1 \cup V_2} \wotimes^{\mathbb{L}}_{A}B\to (A_{V_1} \wotimes_{A}B)\times (A_{V_2} \wotimes_{A}B) \to A_{V_1 \cap V_2}\wotimes_{A}B
\]
and so applying Lemma \ref{lem:2of3discrete} we are finished.
\ \hfill $\Box$

Thus, now we can state our first result of this section.

\begin{thm}\label{thm:LocToHoEp}
Suppose that the morphism $f: \A \to \B$ of $k$-dagger affinoid algebras is a $k$-dagger affinoid localization of a dagger affinoid algebra $A$. Then, for any morphism of $k$-dagger affinoid algebras $A\to C$, the natural morphism $C\wotimes^{\mathbb{L}}_{A}B  \to C\wotimes_{A}\B$ is an isomorphism. Therefore by taking $C=  B$ we see that $f$ is a homotopy epimorphism. 
\end{thm}
{\bf Proof.}
By the main corollary of the Gerritzen-Grauert theorem for $k$-dagger affinoid spaces mentioned in Theorem \ref{thm:MainResults} we may assume that $B=A_{V}$ where $V=V_1 \cup V_2 \cup \cdots \cup V_n$ and the $V_i$ are rational $k$-dagger domains. Now by Remark \ref{rem:usual} (3) we know that any intersection of rational domains is rational. Then by Lemma \ref{lem:CapCup} we are reduced to the case that $f:A \to B$ is a $k$-dagger rational localization of $B$.

\ \hfill $\Box$

Our next goal is to prove the other implication, \ie that any homotopy epimorphism of dagger affinoid algebras is a dagger affinoid localization.

\begin{lem} \label{lemma_splitting}
Let $A, B, C$ be $k$-dagger affinoid algebras considered as objects in $\ttComm(\ttCBorn_k)$ and let $f : \A \to C$ be morphism in $\ttComm(\ttCBorn_k)$ which is a strict epimorphism in $\ttCBorn_k$ such that the post-composition of $f$ with some homotopy epimorphism $g: C \to \B$  in $\ttComm(\ttCBorn_k)$ is a homotopy epimorphism $h : \A \to \B$  in $\ttComm(\ttCBorn_k)$. Then there is a $k$-dagger affinoid algebra $\A'$ and an isomorphism $A \cong C \times \A'$ such that the projection to $C$ corresponds to $f$ under this isomorphism.
\end{lem}

{\bf Proof.} Consider the composition $A \stackrel{f} \longrightarrow C \stackrel{g}\longrightarrow B$. It induces morphisms 

\[D^{\leq 0}(B) \stackrel{}{\to} D^{\leq 0}(C) \stackrel{}{\to} D^{\leq 0}(A),
\]
and since the composition and first morphism are fully faithful (using that $g$ and $g\circ f$ are homotopy epimorphims) then the second must be as well. Thus also $f$ is a homotopy epimorphism and so, by Lemma \ref{lem_HomotopyMon}, we deduce that $C \wotimes^{\Lb}_\A C \cong C$. If $I = \ker (f)$, there is a strict, short exact sequence 
\begin{equation}\label{eqn:sses} 0 \to I \to \A \to C \to 0 \end{equation}
of complete bornological $k$-vector spaces. Applying the functor $V \mapsto V \wotimes^{\Lb}_\A C$ to (\ref{eqn:sses}) we obtain the sequence exact triangle
\[ I \wotimes^{\Lb}_\A C \to C \to C \wotimes^{\Lb}_\A C. \]
But last morphism is an isomorphism, hence $I \wotimes^{\Lb}_\A C = 0$. Now, applying the functor $V \mapsto V \wotimes^{\Lb}_\A I$ to (\ref{eqn:sses}) we get the exact triangle
\[ I \wotimes^{\Lb}_\A I \to I \to C \wotimes^{\Lb}_\A I = 0, \]
which shows that the map $I \wotimes^{\Lb}_\A I \to I$ is an isomorphism. In particular this implies that $I \wotimes_\A I \to I$ is surjective and so also $I \otimes_\A I \to I$, because $I$ is a finitely generated $\A$-modules and hence $I \wotimes_\A I \cong I \otimes_\A I$.
Therefore, $I^{2}=I$ and so there exists an idempotent element $e \in A$ such that $e A = I$ that induces of $I$ a structure of a $k$-dagger affinoid algebra, which we denote 
by $A'=A/(1-e)A$. This induces a splitting of the exact sequence 
\[0 \to I \to A \stackrel{f}\to C \to 0, \]
thus the map $(e,f):A \to A' \times C$ is an isomorphism.
\ \hfill $\Box$

\begin{lem} \label{lemma_last}
Let $\A, \B$ be $k$-dagger affinoid algebras, let $f: A \to B$ be a morphism in the category of $k$-dagger affinoid algebras and let $\{ \mM(A_{V_i})= V_i , i \in I\}$ be a finite dagger affinoid covering of $\mM(A)$. Suppose that for any $V_i$ the morphism $\A_{V_i} \to \A_{V_i} \wotimes_\A \B$ is a $k$-dagger affinoid localization then also $f$ is a $k$-dagger affinoid localization.\end{lem}

{\bf Proof.}

Since dagger affinoid immersion are stable by pullbacks, then $B \to \A_{V_i} \wotimes_\A \B$ is a dagger affinoid immersion and hence the family $\{ \mM(\A_{V_i} \wotimes_\A \B) \}$ is a dagger affinoid covering of $\mM(B)$. Moreover, since $\A_{V_i} \to \A_{V_i} \wotimes_\A \B$ is a dagger affinoid immersion, then also the composition $A \to \A_{V_i} \to \A_{V_i} \wotimes_\A \B$ and so $U \subset \mM(A)$, the image of $\mM(B)$, can be covered by a finite covering of dagger affinoid subdomains of $A$. By Tate acyclicity theorem we have that
\[ \B \cong \A_U = \ker(\prod_{i} \A_{U \cap V_i} \to \prod_{i,j} \A_{U \cap V_i \cap V_j}), \]
which by Proposition 6.1.27 of \cite{Bam} implies that $f$ is a $k$-dagger affinoid localization.
\ \hfill $\Box$

\begin{thm}\label{thm:HoEpToEmb}
Let $\A, \B$ be $k$-dagger affinoid algebras and let $f : \A \to \B$ be a morphism in the category of $k$-dagger affinoid algebras. Assume that $f$ is a homotopy epimorphism 
when considered in the category $\ttComm(\ttCBorn_k)$. Then the morphism $\mM(\B) \to \mM(\A)$ corresponding to $f$ is a $k$-dagger affinoid domain immersion.\end{thm}

{\bf Proof.}
By applying the Gerritzen-Grauert theorem for $k$-dagger affinoid spaces \cite{Bam}, to the morphism $f: \A \to \B$ with the weaker hypothesis that $f$ is an epimorphism we obtain $k$-dagger affinoid algebras rational localization $\A \to \A_{V_i}$ associated to a finite rational covering $\{V_i\}$ of $\mM(\A)$ for which the canonical morphism $\A_{V_i} \to \A_{V_i} \wotimes_\A \B$ corresponds geometrically to a Runge immersion. This means that the map $\A_{V_i} \to \A_{V_i} \wotimes_\A \B$ factorizes as 
\[\A_{V_i} \twoheadrightarrow C_i \hookrightarrow \B \wotimes_\A \A_{V_i}\]
where the first map is the quotient by an ideal and the second is a Weierstrass localization. Now, we use the fact that $f$ is a homotopy epimorphism (and not just an epimorphism), \ie that
\[ \B \wotimes^\Lb_\A \B \cong \B . \]
Since homotopy epimorphism are closed by derived base change we get that
\begin{equation}\label{eqn:something}(\B \wotimes^\Lb_\A \A_{V_i})\wotimes^\Lb_{\A_{V_i}}(\B \wotimes^\Lb_\A \A_{V_i}) \cong \B \wotimes^\Lb_\A \A_{V_i}. \end{equation}
By Lemma \ref{lemma_subdomain} we get $\B \wotimes^\Lb_\A \A_{V_i} \cong \B \wotimes_\A \A_{V_i}$ and hence using equation (\ref{eqn:something}) we conclude that
\[(\B \wotimes_\A \A_{V_i})\wotimes^\Lb_{\A_{V_i}}(\B \wotimes_\A \A_{V_i}) \cong \B \wotimes_\A \A_{V_i}. \]
Thus, the morphism $\A_{V_i} \to \B \wotimes_\A \A_{V_i}$ is a homotopy epimorphism for any $i$.

We already showed in Lemma \ref{lem:WeierLau} that $k$-dagger Weierstrass localizations are homotopy epimorphisms, hence we can apply Lemma \ref{lemma_splitting} to the sequences
\[ \A_{V_i} \twoheadrightarrow C_i \hookrightarrow \B \wotimes_\A \A_{V_i} \]
obtaining dagger affinoid algebras $A'_i$ such that $\A_{V_i} \cong C_i \times \A'_i$  (and so $\mM(\A_{V_i}) \cong \mM(C_i) \coprod (\A'_i)$). Therefore, the map $\A_{V_i} \twoheadrightarrow C_i$ is in fact a $k$-dagger affinoid localization and so also $\A_{V_i} \to \B \wotimes_\A \A_{V_i}$ is a $k$-dagger affinoid localization because is a composition of two $k$-dagger affinoid localizations. We end the proof noticing that the map $f$ and the maps $\A_{V_i} \to \B \wotimes_\A \A_{V_i}$ we discussed satisfies the conditions of Lemma \ref{lemma_last}, so $f$ is a $k$-dagger affinoid localization. Therefore, the morphism $\mM(\B) \to \mM(\A)$ corresponding to $f$ is a $k$-dagger affinoid domain immersion.
\ \hfill $\Box$
\begin{rem}
When $k$ is non-Archimedean we can reduce (using reindexing) the dagger affinoid case to the standard non-Archimedean affinoid one dealt with in \cite{BeKr}. This would give us a quick proof Theorems \ref{thm:LocToHoEp} and \ref{thm:HoEpToEmb} in that case. But, since this reduction step to the affinoid picture cannot be worked out when $k$ is Archimedean, this forces us to reproduce the arguments in the new context. The lack of a reasonable explicit theory of affinoid localizations in the Archimedean context is in fact a major reason to use dagger affinoid geometry and dagger affinoid localizations.
\end{rem}

We conclude this section with saying that we expect that the abstractly defined Grothendieck topology described in subsection \ref{Topologies} restricts to the standard weak G-topology on the category of $k$-dagger affinoid spaces. The covers are clearly invariant by base change and composition and the equivalence of the surjectivity condition on covers with the conservativity requirement is similar to work done in \cite{BeKr}.

\section{A sketch of a theory of dagger affinoid spaces over $\Zb$}\label{sketch}

In this section we sketch some ideas for generalizing the theory of dagger affinoid spaces over a general Banach base ring. Our main aim is to show how overconvergent analytic functions form a natural framework where to work out such a theory, showing also that some differences between Archimedean and non-Archimedean base rings automatically disappear when working with an overconvergent approach.

Consider the following classical settings: the algebras
\[ S_\mathbb{R}^n(\rho) = \{ \sum_{I \in \mathbb{Z}^{n}_{\geq 0}} a_I X^I \ \ | \ \  \sum_{I \in \mathbb{Z}^{n}_{\geq 0}}|a_I| \rho^I < \infty , a_I \in \mathbb{R}\} \]
for a polyradius $\rho = (\rho_1, \dots, \rho_n)$. On $S_\mathbb{R}^n(\rho)$ one can consider the spectral semi-norm defined
\[ |f|_\sup = \max_{|\cdot| \in \mM(S_\mathbb{R}^n(\rho))} |f(x)| \]
for any $f \in S_\mathbb{R}^n(\rho)$. We denote by $T_\mathbb{R}^n(\rho)$ the completion of $S_\mathbb{R}^n(\rho)$ with respect to this semi-norm, which is in fact a norm. The algebras $T_\mathbb{R}^n(\rho)$ are classically called \emph{disk algebras} and elements of $T_\mathbb{R}^n(\rho)$ correspond with continuous function on the closed polydisk of polyradius $\rho$ which are analytic in the interior. Then, the following isomorphism of complete bornological algebras is well known
\[ \limind_{\rho' > \rho} S_\mathbb{R}^n(\rho') \cong \limind_{\rho' > \rho} T_\mathbb{R}^n(\rho'). \]
The same holds for $\mathbb{R}$ replaced by any Archimedean complete valuation field. Quite surprisingly, this is not a peculiar feature of the Archimedean context.

Thus, let's consider the analogous algebras over $k$, with $k$ non-Archimedean. We define
\begin{equation}\label{eqn:SKN} S_k^n(\rho) = \{ \sum_{I \in \mathbb{Z}^{n}_{\geq 0}} a_I X^I\ \  | \ \  \sum_{I \in \mathbb{Z}^{n}_{\geq 0}} |a_I| \rho^I < \infty , a_I \in k\} \end{equation}
which is an object of $\ttComm(\ttBan^{A}_{k})$, whose elements can be interpreted as functions on
\[ \{c\in k^{n} \ \ | \ \ |c_{i}| \leq \rho_{i} \ \ 1\leq i \leq n\} \] 
of polyradius $\rho$. We equip $S_k^n(\rho)$ with the norm 
\begin{equation}\label{eqn:SumNorm}\|\sum_{I \in \mathbb{Z}^{n}_{\geq 0}} a_I X^I\| = \sum_{I \in \mathbb{Z}^{n}_{\geq 0}} |a_I|\rho^{I}.
\end{equation} Note that the norm of $S_k^n(\rho)$ is not non-Archimedean and hence the topology induced by the norm on $S_k^n(\rho)$ is not locally convex. Since it is unusual to consider this kind of algebras over non-Archimedean base fields, we spend some words about their basic properties, for which the only reference in literature we know is the book \cite{GR2}, which has not an english translation. 

\begin{prop}
The sets we defined in equation (\ref{eqn:SKN}) $S_k^n(\rho)$ form well defined subalgebras of $k[[X_1, \ldots, X_n]]$ which are complete when equipped with the norm from Equation (\ref{eqn:SumNorm}).
\end{prop}
{\bf Proof.}
The fact that $S_k^n(\rho)$ are closed under addition and multiplication follows readily by the triangle inequality, which also shows that $\|f g\| \le \|f\|\|g\|$ for any $f,g \in S_k^n(\rho)$. Hence the only non-trivial fact to check is the completeness with respect to the norm.

Let's consider a Cauchy sequence $\{ f_j = \sum_{I \in \mathbb{Z}^{n}_{\geq 0}} a_{j, I} X^I  \}$ in $S_k^n(\rho)$. Then for each $I$ we can consider the Cauchy sequences $a_{j, I}$ of coefficients in $k$ which admit a limit
\[ \lim_{j \to \infty} a_{j, I} = a_I \]
since $k$ is complete. Then, the fact that the power-series
\[ f = \sum_{I \in \mathbb{Z}^{n}_{\geq 0}} a_I X^I \]
is a limit of the Cauchy sequence $\{ f_j \}$ can be proved by standard arguments noticing that that 
\[ \lim_{j \to \infty} \| f - f_j\| = 0 \]
and that $\|\cdot\|$ is a norm.
\ \hfill $\Box$

Let now $k$ be a non-trivially valued non-Archimedean field. Consider  the Tate algebra
\[
T_k^n(\rho) =k \lt \rho_{1}^{-1}x_{1}, \dots, \rho_{n}^{-1} x_{n}\gt = \{ \sum_{I \in \mathbb{Z}^{n}_{\geq 0}} a_{I}x^{I} \in k[[x_1, \dots, x_n]] \ \ | \ \ \lim_{|I| \to \infty} |a_{I}| \rho^{I} = 0 \}
\]
where $|I|= i_{1} + \cdots +i_{n}$ and we equip the right hand side with the norm 
\[\|\sum_{I \in \mathbb{Z}^{n}_{\geq 0}} a_{I} x^{I}\|_{\rho}=\sup_{I\in \mathbb{Z}^{n}_{\geq 0}}\{|a_{I}|\rho^{I}\}.
\]

One can think of $T_k^n(\rho)$ as the completion of $S_k^n(\rho)$ with respect to this norm, which coincides with the spectral norm of $S_k^n(\rho)$. In the analogous situation over Archimedean base fields the completion of $S_k^n(\rho)$ is the algebra of analytic function on the polycylinder of polyradius $\rho$ which are continuous on the boundary. So, the impression is that the algebra $S_k^n(\rho)$ is the more fundamental, from which we can recover all the algebras of analytic geometry and the only one which has a perfectly uniform description over any valued field. Moreover, its definition also naturally generalize over any Banach ring.

Despite the fact that $S_k^n(\rho)$ and $T_k^n(\rho)$ have different properties (one is locally convex and the other is not), the following interesting result holds.

\begin{thm} \label{thm:non_arch_summation_dagger}
Let $k$ be a non-Archimedean valued field and $r$ any polyradius then
\[ \limind_{\rho > r} S_k^n(\rho) \cong \limind_{\rho > r} T_k^n(\rho) \]
as objects in $\ttComm(\ttCBorn_{k})$ where this colimit is taken in $\ttComm(\ttCBorn_{k})$ .
\end{thm}

{\bf Proof.}

It is enough to show that the two systems are final in each other. Consider any polyradius $\rho < \rho'$ (\ie each component of $\rho$ is strictly less than the corresponding component of $\rho'$) we need to show that
\[  S_k^n(\rho') \subset T_k^n(\rho) \]
and the other inclusion
\[ T_k^n(\rho') \subset S_k^n(\rho). \]
Let $f = \sum_{I \in \mathbb{Z}^{n}_{\geq 0}} a_I X^I \in S_k^n(\rho')$ then
\[ \sum_{I \in \mathbb{Z}^{n}_{\geq 0}} |a_I| (\rho')^I < \infty \then |a_I|(\rho')^I \to 0 \then f \in T_k^n(\rho') \then f \in T_k^n(\rho). \]
On the other hand, let $f = \sum_{I \in \mathbb{Z}^{n}_{\geq 0}} a_I X^I \in T_k^n(\rho')$, then  $|a_I|(\rho')^I \to 0$, hence then series
\[ \sum_{I \in \mathbb{Z}^{n}_{\geq 0}} |a_I| r^I \]
converges for any $r < \rho'$. 

This shows that there is a bijection between the two systems, which readily implies an isomorphism of bornological vector spaces because the map we considered above are bounded. More precisely, let's consider the restriction map
\[ S_k^n(\rho') \to T_k^n(\rho). \]
This map is bounded since factor through $S_k^n(\rho') \to T_k^n(\rho')$ and
\[ \max_{I \in \mathbb{Z}^{n}_{\geq 0}} |a_I|(\rho')^I  \le \sum_{I \in \mathbb{Z}^{n}_{\geq 0}}  |a_I|(\rho')^I \]
and the canonical map $T_k^n(\rho') \to T_k^n(\rho)$ is obviously bounded. Finally the map
\[ T_k^n(\rho') \to S_k^n(\rho) \]
is bounded because
\[ \sum_{I \in \mathbb{Z}^{n}_{\geq 0}}  |a_I| \rho^I \le \max_{1 \le i \le n} \left ( \frac{\rho_i'}{\rho_i' - \rho_i} \right ) \max_{I \in \mathbb{Z}^{n}_{\geq 0}}   |a_I|(\rho')^I  \]
showing that the bijections are isomorphisms of bornological algebras.
\ \hfill $\Box$

\begin{rem}
We remark that the boundedness of the map $T_k^n(\rho') \to S_k^n(\rho)$, if $k$ is Archimedean and $T_k^n(\rho')$ is, as in the first part of this section, the disk algebra of the unital polydisk, is a non-trivial consequence of Cauchy integration formula. Instead, for $k$ non-Archimedean this result is very easy, as explained above, and this is due to the fact that the Shilov boundary of the polydisk is made of a single point.
\end{rem}

Another, property of these kind of systems is that the previous result holds also if the base field is trivially valued, provided that it is interpreted in the language of Ind-objects. In fact in the proof of Theorem \ref{thm:non_arch_summation_dagger} there is nothing depending on the fact that the base field is non-trivially valued and so the proof works for any non-Archimedean complete valued field. 



Let's consider the analogous statements and algebras for rings. So, let's define
\[ S_R^n(\rho) = \{ \sum_{I \in \mathbb{Z}^{n}_{\geq 0}} a_I X^I  \ \ | \ \ \sum_{I \in \mathbb{Z}^{n}_{\geq 0}}|a_I| \rho^I < \infty , a_I \in R \} \]
equipped with the norm $\|\sum_{I \in \mathbb{Z}^{n}_{\geq 0}} a_I X^I\| = \sum_{I \in \mathbb{Z}^{n}_{\geq 0}} |a_I| \rho^I$. $S_R^n(\rho)$ is in fact a Banach ring (and in particular a Banach $R$-algebra), so by Theorem 1.2.1 of \cite{Ber1990} its spectrum $\mM(S_R^n)$ is a non-empty, compact, Hausdorff topological space. The elements of $\mM(S_R^n)$ (or in general of any Banach ring) are bounded multiplicative semi-norms and $\mM(S_R^n)$ is equipped with the weakest topology for which the maps $|\cdot| \mapsto |f|$ are continuous for every $f \in R$.

\begin{defn}
Let $R$ be a Banach ring, then the \emph{spectral semi-norm} of $R$ is defined
\[ |f|_\sup = \max_{|\cdot| \in \mM(R)} |f|, \]
for any $f \in R$.
\end{defn}

One can check that 
\[ |f|_\sup = \lim_{n \to \infty} \sqrt[n]{|f|^n} = \inf \sqrt[n]{|f|^n}. \]
\begin{defn}
We denote the separated completion of $S_R^n(\rho)$ with respect to the spectral semi-norm by $T_R^n(\rho.)$ 
\end{defn}
Note that $T_R^n(\rho)$ is a Banach $R$-algebra. In order to define the over-convergent functions on the unit disk, we can consider the objects
\[ ``\limind"_{\rho > 1} S_R^n(\rho) \]
and
\[ ``\limind"_{\rho > 1} T_R^n(\rho) \]
of $\ttComm(\ttInd(\ttBan_R)$. We remark that all the maps that defines 
these objects are injective, hence the underlying $\ttInd(\ttBan_R)$-module is an essentially monomorphic object, as defined in Definition \ref{defn:essential_monomorphic}, of $\ttInd(\ttBan_R)$. Thus, since this category is concrete (Corollary \ref{prop:essential_concrete}) it is meaningful to associate to $``\limind"_{\rho > 1} S_R^n(\rho)$ and $``\limind"_{\rho > 1} T_R^n(\rho)$ their underlying rings and to think of these inductive systems as a structure over this ring given by a family of norms, like a bornology. In fact we can even give the following definition.

\begin{defn}
Let $R$ be a Banach ring, we define the following categories:
\begin{itemize}
\item $\ttBorn_R$, the category of \emph{bornological modules of convex type} over $R$, to be the sub-category of essential monomorphic objects of $\ttInd(\ttSNrm_R)$;
\item $\ttSBorn_R$, the category of \emph{separated bornological modules} over $R$, to be the sub-category of essential monomorphic objects of $\ttInd(\ttNrm_R)$;
\item $\ttCBorn_R$, the category of \emph{complete bornological modules} over $R$, to be the sub-category of essential monomorphic objects of $\ttInd(\ttBan_R)$.
\item when $R$ is non-Archimedean we can take $\ttCBorn^{nA}_R$, the category of \emph{complete non-
Archimedean bornological modules} over $R$, to be the sub-category of essential monomorphic objects of $\ttInd(\ttBan^{nA}_R)$.
\end{itemize}
\end{defn}
\begin{rem}
In the case that $R$ is a complete valued field, this definition is equivalent to the categories  of Archimedean or non-Archimedean bornological vector spaces which we used for a non-trivially valued complete field. In general, the category $\ttCBorn_{R}$ is an elementary quasi-abelian closed symmetric monoidal category, has all limits and colimits and enough flat projectives so it can be used for relative algebraic geometry.
\end{rem}
 We are particulary  interested in studying the case when $R = \Zb$. We remark that for $\rho \ge 1$ there are algebraic isomorphisms 
\[ S_\Zb^n(\rho) \cong \Zb[X_1, \ldots X_n] \]
\[ T_\Zb^n(\rho) \cong \Zb[X_1, \ldots X_n]. \]
So, one might be lead to think that objects like $``\limind_{\rho > 1}" S_\Zb^n(\rho)$ and $ ``\limind_{\rho > 1}" T_\Zb^n(\rho)$
are not worthy of study. However, this is not true and we will see that the family of norms that defines these inductive systems encodes important analytic information of the global unital disk.

Recall that to each $x \in \mM(\Zb)$ one can associate the residue field $\Hc(x)$ by extending the semi-norm associated to $x$, denoted $|\cdot|_x$, to a valuation on $\Zb /\ker(|\cdot|_x)$ and then complete it. Hence, each $\Hc(x)$ is a complete valued field which comes canonically associated with a bounded morphism $(\Zb, |\cdot|_\infty) \to \Hc(x)$.

We will need the following lemma, which describes the Shilov boundary of the unital polydisk of the affine global spaces over $\Zb$ in the sense of Poineau. We remark that in \cite{Poi2} there are similar descriptions of Shilov boundaries of compact subsets of $\Ab_\Zb^{n, \text{an}}$, but the study of Shilov boundary done in next lemma is missing. In this lemma, we show that the morphism 
\begin{equation}\label{eqn:themor}\mM(S_\Zb^n(\rho) \wotimes_\Zb \Rb) \to \mM(S_\Zb^n(\rho))
\end{equation}
 induced by $S_\Zb^n(\rho) \to S _\Zb^n(\rho) \wotimes_\Zb \Rb$ identifies Shilov boundaries. Here. $\Zb$ and $\Rb$ are treated as Archimedean Banach rings with the norm being the standard absolute values.

\begin{lem}
The space $\mM(S_\Zb^n(\rho))$ is the polydisk of polyradius $\rho$ in the global affine space of Poineau and the Shilov boundary of $\mM(S_\Zb^n(\rho))$ agrees with the Shilov boundary of $\mM(S_\Zb^n(\rho) \wotimes_\Zb \Rb)$ via (\ref{eqn:themor}).
\end{lem}

{\bf Proof.}

One can check that set-theoretically 
\[ \mM(S_\Zb^n(\rho)) = \coprod_{x \in \mM(\Zb)} \mM(S_\Zb^n(\rho) \wotimes_{\mathbb{Z}} \Hc(x)) \]
and $S_\Zb^n(\rho) \wotimes_{\mathbb{Z}} \Hc(x) \cong S_{\Hc(x)}^n(\rho)$.
Then, since points of $\mM(S_\Zb^n(\rho)))$ are bounded multiplicative semi-norms it is clear that for any Archimedean point $x \in \mM(\Zb)$ whose complete residue field is not isometrically isomorphic to $(\Rb, |\cdot|_\infty)$ the sets $\mM(S_\Zb^n(\rho) \wotimes_{\mathbb{Z}} \Hc(x))$ cannot contain any point of the Shilov boundary of $\mM(S_\Zb^n(\rho))$. And this is true also for all point with non-Archimedean residue fields because in that case we have the easy estimate
\[ \max_{z \in \mM(S_\Zb^n(\rho) \wotimes_{\mathbb{Z}} \Hc(x))} |f(z)| \le 1  \] 
for any non-Archimedean $x \in \mM(\Zb)$. Hence the only possibility is that the Shilov boundary of $\mM(S_\Zb^n(\rho))$ and $\mM(S_\Zb^n(\rho) \wotimes_{\Zb} \Rb)$ are identified.

\ \hfill $\Box$

In the next result we see that different possible definitions of dagger algebras agree, using the usual absolute value on $\Zb$.

\begin{prop}\label{prop:SameZ}
In $\ttComm(\ttInd(\ttBan_{\Zb}))$ or $\ttComm(\ttCBorn_{\Zb})$ there is an isomorphism
\[ ``\limind_{\rho > r}" S_\Zb^n(\rho) \cong ``\limind_{\rho > r}" T_\Zb^n(\rho), \]
for any polyradius $r$.
\end{prop}
{\bf Proof.}

There is a canonical morphism
\[ \phi = ``\limind"_{\rho > r} \phi_\rho: ``\limind_{\rho > r}" S_\Zb^n(\rho) \to ``\limind_{\rho > r}" T_\Zb^n(\rho) \]
where $\phi_{\rho}:S_\Zb^n(\rho)\to T_\Zb^n(\rho) $ is obtained by completing each element of the system on the left side with respect to the spectral norm. This map is clearly a monomorphism, so it is enough to show that this map has an inverse. As for Theorem \ref{thm:non_arch_summation_dagger} it is enough to show that for any $\rho < \rho'$ the restriction map
\[ T_\Zb^n(\rho') \to S_\Zb^n(\rho) \]
is well defined and bounded. By previous lemma this is reduced to the known case over $\Cb$ because (as a consequence of the Cauchy formula in several variables) one can show that the restriction morphism
\[ \psi_{\rho'}: T_\Cb^n(\rho') \to S_\Cb^n(\rho) \]
is bounded and it restricts to a bounded morphism $T_\Zb^n(\rho') \to S_\Zb^n(\rho)$, which is injective. A morphism 
\[\psi= \limind \psi_{\rho}: ``\limind"_{\rho > r} T_\Zb^n(\rho) \to ``\limind"_{\rho > r} S_\Zb^n(\rho)
\] can be defined by properly re-indexing both systems, considering two sequences of polyradii $\{\rho_n\}$ and $\{\rho'_n\}$ both converging to $r$ and with $\rho_n' > \rho_n$ for every $n$. It is clear that such a choice is always possibile.

To conclude that $\phi$ and $\psi$ are inverse of each other is enough to check that they induces a bijection of sets $U(``\limind"_{\rho > r} S_\Zb^n(\rho)) \cong U(``\limind"_{\rho > r} T_\Zb^n(\rho))$, where $U$ is the concretization functor defined as before Corollary \ref{prop:essential_concrete}, noticing that both systems are monomorphic. This bijection is readily deduced by the inclusion of sets
\[ S_\Zb^n(\rho') \subset T_\Zb^n(\rho') \subset S_\Zb^n(\rho) \]
and 
\[ T_\Zb^n(\rho') \subset S_\Zb^n(\rho) \subset T_\Zb^n(\rho) \]
for any $\rho < \rho'$.

\ \hfill $\Box$

\begin{rem}
The previous result can also be done in more general settings using Proposition 2.1.3 of \cite{Poi2} and Corollary 2.8 of \cite{Poi3}. We preferred to give an easy proof of the special case of disks over $\Zb$ to provide the reader with an explicit key example to train intuition. Moreover, even if with Poineau's results one can deduce the previous one, the Ind-objects point of view is not present in \cite{Poi2} and \cite{Poi3}. We think that this interpretation can give a useful geometric insight toward the understanding of global analytic spaces.
\end{rem}

Hence also in the case of $\Zb$ the choice of considering the overconvergent power-series remove the embarrassment of the choice of norm to put on the polynomial to get the right ring to study. Both natural choices leads to the same objects that we can denote, unambiguously, with the symbol
\[ \Zb \lt \rho_1^{-1} X_1, \ldots, \rho_n^{-1} X_n \gt^\dagger. \]

%
%
%
%

The previous discussion was motivated by the following problem. The category of global analytic spaces over $\mathbb{Z}$ introduced by Poineau in \cite{Poi2} allows one to work only with spaces which are ``without boundary", \ie the building blocks are not some sort of affinoid spaces but are spaces defined by subsets of opens of affine spaces identified by the zeros of a finite number of analytic functions. This is a definition of analytic spaces in the spirit of complex analytic geometry. 
Instead one can ask to build a theory which is some sort of generalization of affinoid spaces by considering closed polydisks and spaces defined by zeros of analytic functions on them. There are several difficulties involved in constructing such a theory and we will not address all of them in this work. Instead, we show that what has been done up to here can be helpful towards this goal.

\begin{thm} \label{thm:CompareZR}
Let $\Zb \lt \rho_1^{-1} X_1, \ldots, \rho_n^{-1} X_n \gt^\dagger$ then
\[ \Zb \lt \rho_1^{-1} X_1, \ldots, \rho_n^{-1} X_n \gt^\dagger \wotimes_\Zb \Rb \cong \Rb \lt \rho_1^{-1} X_1, \ldots, \rho_n^{-1} X_n \gt^\dagger \]
and for any prime $p$ we have,
\[ \Zb \lt \rho_1^{-1} X_1, \ldots, \rho_n^{-1} X_n \gt^\dagger \wotimes_\Zb \Qb_p \cong \Qb_p \lt \rho_1^{-1} X_1, \ldots, \rho_n^{-1} X_n \gt^\dagger \]
\end{thm}

{\bf Proof.}

Since for objects in $\ttInd(\ttBan_{\Zb})$ the tensor product is calculated term by term we get that
\[ \Zb \lt \rho_1^{-1} X_1, \ldots, \rho_n^{-1} X_n \gt^\dagger \wotimes_\Zb \Rb = (\limind_{r > \rho} S_\Zb^n(r))\wotimes_\Zb \Rb \cong \limind_{r > \rho} S_\Rb^n(r) = \Rb \lt \rho_1^{-1} X_1, \ldots, \rho_n^{-1} X_n \gt^\dagger\]
and
\begin{equation}\begin{split} \Zb \lt \rho_1^{-1} X_1, \ldots, \rho_n^{-1} X_n \gt^\dagger \wotimes_\Zb \Qb_p & =  (\limind_{r > \rho} S_\Zb^n(r))\wotimes_\Zb \Qb_p \cong \limind_{r > \rho} S_{\Qb_p}^n(r) \\ & \cong \limind_{r > \rho} T_{\Qb_p}^n(r)= \Qb_p \lt \rho_1^{-1} X_1, \ldots, \rho_n^{-1} X_n \gt^\dagger. \end{split}
\end{equation}
Where we have used Theorem \ref{thm:non_arch_summation_dagger} to show the second isomorphism.
\ \hfill $\Box$

\begin{rem}
Theorem \ref{thm:CompareZR} is not true for non-overconvergent analytic functions. This theorem could also have been proven using Proposition \ref{prop:SameZ}.
\end{rem}

So it seems natural to develop a theory of dagger analytic spaces over any Banach ring (non-Archimedean or not) $R$ using as building block algebras quotients of the algebras
\[ R \lt \rho_1^{-1} X_1, \ldots, \rho_n^{-1} X_n \gt^\dagger = {``\limind"}_{r > \rho} S_R^n(r) \]
for any polyradius $\rho$ by ideals. Or one can restrict to considering only \emph{strict spaces}, in the sense of admitting only $\rho = 1$ as in classical rigid geometry. We expect many of the results of this article to carry over to that case. The category of spaces we obtain will be pretty similar to the one recently obtained in \cite{Pau2} by Paugam.  In the appendix, we define non-Archimidification, a way to produce an interesting functor 
\[\ttComm(\ttInd(\ttBan^{A}_{R})) \to \ttComm(\ttInd(\ttBan^{nA}_{R}))
\]
when $R$ is a non-Archimedean Banach ring. This gives an interesting way of producing non-Archimedean geometry in the standard sense out of the newer type of geometry we are proposing. This morphism sends the standard affine models (quotients of $S^{n}_{R}(\rho)$) in our new type of geometry, to the standard affine models (quotients of $T^{n}_{R}(\rho)$) in standard affinoid theory.


\section{Appendix: Non-Archimidification functors}

\begin{lem}
Given categories $\ttC_{1}$ and $\ttC_{2}$, such that $\ttC_{1}$ has a generator, small colimits and finite limits, and such that small filtrant inductive limits in  $\ttC_{1}$ are stable by base change then any functor
\[\iota:\ttC_{1} \to \ttC_{2}
\]
has a left adjoint.
\end{lem}
{\bf Proof.}
We can define a functor 
\[\ttC^{op}_{1} \to \ttSets
\]
by 
\[W \mapsto \Hom_{\ttC_2}(V,\iota(W)).
\]
By Theorem 5.3.9 of \cite{KS} such functor is representable, that fact clearly provides the desired left adjoint $\pi$. 

\hfill $\Box$

\begin{thm}\label{thm:LeftAdjContr}
For any non-Archimedean valuation ring $R$, the functor
\[\iota: \ttBan^{nA, \leq 1}_{R} \to \ttBan^{A, \leq 1}_{R}
\]
has a left adjoint $\pi$ which respects the monoidal structures. We also have $\pi(P^{A}(V))= P^{nA}(V).$
\end{thm}
{\bf Proof.}
In this proof we will use the notation $\coprod$ to denote the copropduct in $\ttBan^{nA, \leq 1}_{R}$  (we will not write it as $\coprod^{\le 1}$). In this setting notice that for any family $\{ V_{i} \}_{i \in I}$ of objects of $\ttBan^{nA, \leq 1}_{R}$ we have a canonical morphism in $\ttBan^{A, \leq 1}_{R}$ 
\[\coprod_{i \in I} \iota(V_{i}) \to  \iota(\coprod_{i \in I} V_{i}).
\]
So for $V \in \ttBan^{nA, \leq 1}_{R}$ and $\kappa_V^{nA}: P^{nA}(V) \to V$ we have
\[\Hom_{\ttBan^{nA, \leq 1}_{R}}(\coprod_{w \in (\ker \kappa_{V}^{nA})^{\times}}R_{\|w\|} , \coprod_{v\in V^{\times}}R_{\|v\|})=\underset{w\in (\ker \kappa_{V}^{nA})^{\times}}{{{\times}}} \coprod_{v \in V^{\times}}R_{\|w\|^{-1}\|v\|}
\]
where the coproduct on the right is in the non-Archimedean contracting category and similarly for $\kappa_V^{A}: P^{A}(V) \to V$, 

\[\Hom_{\ttBan^{A, \leq 1}_{R}}(\coprod_{w \in (\ker \kappa_{V}^A)^{\times}}R_{\|w\|} , \coprod_{v\in V^{\times}}R_{\|v\|})=\underset{w\in (\ker \kappa_{V}^A)^{\times}}{{{\times}}} \coprod_{v \in V^{\times}}R_{\|w\|^{-1}\|v\|}
\]

where the coproduct on the right is in the Archimedean contracting category.
So we get a morphism
\begin{equation}\label{equation:getmap}
\begin{split}
\Hom_{\ttBan^{A, \leq 1}}(\coprod_{w\in \ker \kappa_{V}^{\times}}R_{\|w\|} , \coprod_{v\in V^{\times}}R_{\|v\|}) &  \\  =\underset{w\in \ker \kappa_{V}^{\times}}{{{\times}}} \coprod_{v\in V^{\times}}\iota(R_{\|w\|^{-1}\|v\|})\to  &\underset{w\in \ker \kappa_{V}^{\times}}{{{\times}}} \iota(\coprod_{v\in V^{\times}}R_{\|w\|^{-1}\|v\|}) \\
& = \Hom_{\ttBan^{nA, \leq 1}}(\coprod_{w\in \ker \kappa_{V}^{\times}}R_{\|w\|} , \coprod_{v\in V^{\times}}R_{\|v\|})
\end{split}
\end{equation}
Then by taking the image of 
\[
\coprod_{w\in \ker \kappa_{V}^{\times}}R_{\|w\|} \stackrel{\kappa_{\ker(\kappa^{A}_{V})}}\longrightarrow \ker(\kappa^{A}_{V}) \hookrightarrow \coprod_{v\in V^{\times}}R_{\|v\|}
\]
 under the morphism from Equation (\ref{equation:getmap}) we get an element of 
 \[ \Hom_{\ttBan_{R}^{nA, \leq 1}}(\coprod_{w\in \ker \kappa_{V}^{\times}}R_{\|w\|} , \coprod_{v\in V^{\times}}R_{\|v\|})
 \] which we can use to define a functor 


\[\pi: \ttBan^{A, \leq 1}_{R} \to \ttBan^{nA, \leq 1}_{R}
\]
by 
\[\pi(V) =  \coker[\coprod_{w\in \ker \kappa_{V}^{\times}}R_{\|w\|}\to \coprod_{v\in V^{\times}} R_{\|v\|}]
\]
where coproducts and the cokernel are taken in $\ttBan^{nA, \leq 1}_{R}.$ We now show that $\pi$ is a left adjoint to the inclusion functor $\iota$. In fact, thanks to Lemma \ref{lem:boundedNorm} we have
\begin{equation}
\begin{split}\Hom_{\ttBan^{nA}_{R}}(\pi(V),W)^{\leq r} & = \ker[\prod_{v\in V^{\times}}\Hom_{\ttBan^{nA}_{R}}(R_{\|v\|},W)^{\leq r} \to \prod_{w\in \ker \kappa_V - \{0\}}\Hom_{\ttBan^{nA}_{R}}(R_{\|w\|},W)^{\leq r}] \\
& =  \ker[\prod_{v\in V^{\times}}\Hom_{\ttBan^{A}_{R}}(R_{\|v\|},\iota (W))^{\leq r} \to \prod_{w\in \ker \kappa_V - \{0\}}\Hom_{\ttBan^{A}_{R}}(R_{\|w\|},\iota(W))^{\leq r}]  \\
& = \Hom_{\ttBan^{A}_{R}}(\coker[\coprod_{w\in \ker \kappa_{V}^{\times}}R_{\|w\|} \to  \coprod_{v\in V^{\times}}R_{\|v\|}],\iota(W))^{\leq r} \\
 &= \Hom_{\ttBan^{A}_{R}}(V,\iota(W))^{\leq r}
\end{split}
\end{equation}
Therefore, we can also conclude that for any $r>0$ we have 
\begin{equation}\label{eqn:alsor}\Hom_{\ttBan^{nA}_{R}}(\pi(V),W)^{\leq r} =\Hom_{\ttBan^{A}_{R}}(V,\iota(W))^{\leq r}.
\end{equation}
 and so 
 \begin{equation}\label{eqn:alsoNor}\Hom_{\ttBan^{nA}_{R}}(\pi(V),W) =\Hom_{\ttBan^{A}_{R}}(V,\iota(W)).
\end{equation}
and
 $\pi$ commutes with colimits. It is now easy to see that $\pi$ intertwines the monoidal structures:
 \begin{equation}\label{eqn:piInt}
\begin{split}\Hom_{\ttBan^{nA, \leq 1}_{R}}(\pi(U \ootimes_{\ttBan^{A, \leq 1}_{R}} V),W) & = \Hom_{\ttBan^{A, \leq 1}_{R}}(U \ootimes_{\ttBan^{A, \leq 1}_{R}} V,\iota(W)) \\ & = \Hom_{\ttBan^{A, \leq 1}_{R}}(U ,\uHom_{\ttBan^{A, \leq 1}_{R}}( V,\iota(W))) \\
& =  \Hom_{\ttBan^{A, \leq 1}_{R}}(U ,\iota(\uHom_{\ttBan^{nA, \leq 1}_{R}}(\pi( V),W))) \\
&=  \Hom_{\ttBan^{nA, \leq 1}_{R}}(\pi(U) ,\uHom_{\ttBan^{nA, \leq 1}_{R}}(\pi( V),W)) 
\\ &=  \Hom_{\ttBan^{nA, \leq 1}_{R}}(\pi(U) \ootimes_{\ttBan^{nA, \leq 1}_{R}} \pi( V),W).
\end{split}
\end{equation}

The passage from the second to the third line requires some explanation. Notice first it is true for $U=R_{w}$ because by Equation \ref{eqn:alsor} we have
\begin{equation}
\begin{split}
\Hom^{\leq 1}(R_w, \uHom(V,\iota(W))) =\Hom(V,\iota(W))^{\leq w} & = \Hom(\pi(V),W)^{\leq w}\\ & = \Hom^{\leq 1}(R_w, \uHom(\pi(V),W)).
\end{split}
\end{equation}
For a general object $U$ we can write $U= \coker[P^{A}(\ker(\kappa_{U}) \to P^{A}(U)]$ and use the fact that colimits in the first component of $\Hom$ become limits of the $\Hom$ sets.

By the Yoneda lemma for the opposite category of $\ttBan^{nA, \leq 1}_{R}$ we conclude using \ref{eqn:piInt} for all $W$ that the natural morphism 
\[\pi(U) \ootimes_{\ttBan^{nA, \leq 1}_{R}} \pi( V) \to \pi(U \ootimes_{\ttBan^{A, \leq 1}_{R}} V)
\]

is an isomorphism. Finally,
\begin{equation}\begin{split}\Hom_{\ttBan^{nA, \leq 1}_{R}}(\pi(P^{A}(V)), W) = \Hom_{\ttBan^{A, \leq 1}_{R}}(P^{A}(V), \iota(W)) & = {{\times}}_{v\in V^{\times}}\{w\in W\ \ | \ \ \|w\|\leq \|v\|\} \\
& = \Hom_{\ttBan^{nA, \leq 1}_{R}}(P^{nA}(V), W).
\end{split}
\end{equation}
\ \hfill $\Box$

\begin{lem}\label{lem:same}The canonical morphism $P^{A}(V)\to P^{nA}(V)$ induces for any $W \in \ttBan^{nA}_{R}$ a natural isomorphism 
\[\Hom_{\ttBan^{nA}_{R}}(P^{nA}(V),W)\to \Hom_{\ttBan^{A}_{R}}(P^{A}(V),\iota(W)).\]
\end{lem}


{\bf Proof.}
The morphism preserves norms, so it is enough to check that for every $r>0$ we get a bijection $\Hom_{\ttBan^{nA}_{R}}(P^{nA}(V),W)^{\leq r}\to \Hom_{\ttBan^{A}_{R}}(P^{A}(V),\iota(W))^{\leq r}.$ This is an immediate consequence of Equation \ref{eqn:alsor} and the fact that $\pi(P^{A}(V)) = P^{nA}_{V}$ which was proven as part of Lemma \ref{thm:LeftAdjContr}.

\ \hfill $\Box$
\begin{thm}For any non-Archimidean valuation ring $R$, the natural morphism
\[\iota: \ttInd(\ttBan^{nA}_{R}) \to \ttInd(\ttBan^{A}_{R})
\]
has a left adjoint which respects the monoidal structures.
\end{thm}
{\bf Proof.}
Now, the category $\ttInd(\ttBan^{nA}_{R})$ has a generator (see Definition \ref{def:generator})  in the sense of Kashiwara-Schapira by Lemma \ref{lem:hasgen}. Lemma \ref{lem:exact_direct_limit} implies that small filtered inductive limits in $\ttInd(\ttBan^{nA}_{R})$ are exact in the sense of Definition \ref{def:KSexact}. We need to use this exactness to check that small filtrant inductive limits in  $\ttInd(\ttBan^{nA}_{R})$ are stable by base change. However, this is implied by Lemma 3.3.9 of \cite{KS}.

\ \hfill $\Box$

We now give an explicit construction of the adjoint functor in a way similar to what was done in the non-expanding category. We will use the construction of projectives in $\ttInd(\ttBan^{nA}_{R})$ and $\ttInd(\ttBan^{A}_{R})$ from Lemma \ref{lem:exact_direct_limit} Notice that, 
\[
\mathcal{E} \cong \coker[P^{A}(\ker(\kappa^{A}_{\mathcal{E}}))  \to P^{A}(\mathcal{E})].
\]
Define 
\[\pi(\mathcal{E}) = \coker[P^{nA}(\ker(\kappa^{nA}_{\mathcal{E}}))  \to P^{nA}(\mathcal{E})].
\]
\begin{equation}\label{IndNAa} \begin{split}\uHom_{\ttInd(\ttBan^{nA}_{R})}(\pi(\mathcal{E}), F) & = \ker [\uHom_{\ttInd(\ttBan^{nA}_{R})}(P^{nA}(\mathcal{E}),F) \to \uHom_{\ttInd(\ttBan^{nA}_{R})}(P^{nA}(\ker(\kappa^{nA}_{\mathcal{E}})),F)] \\ & =  \ker [\lim_{J \subset I,  |J|<\infty} \colim_{k \in K} \uHom_{\ttBan^{nA}_{R}}(\coprod_{j\in J}P^{nA}(E_j),F_{k}) \\ & \to \lim_{T \subset L, |T|<\infty} \colim_{p\in P} \uHom_{\ttBan^{nA}_{R}}(\coprod_{t\in T}P^{nA}(\ker(\kappa^{nA}_{\mathcal{E}}))_{t},F_p)] \\
\end{split}
\end{equation}
Using Lemma \ref{lem:same} and because $J$ is finite, we know that
 \begin{equation}\label{IndNAb}
\begin{split}\uHom_{\ttBan^{nA}_{R}}(\coprod_{j\in J}P^{nA}(E_j),F_k)= \prod_{j\in J} \uHom_{\ttBan^{nA}_{R}}(P^{nA}(E_j),F_k) & = \prod_{j\in J} \uHom_{\ttBan^{A}_{R}}(P^{A}(E_j),F_k) \\ & = \uHom_{\ttBan^{A}_{R}}(\coprod_{j\in J}P^{A}(E_j),F_k)
\end{split}
\end{equation}
and similarly, because $T$ is finite,
 \begin{equation}\label{IndNAc}
\begin{split}\uHom_{\ttBan^{nA}_{R}}(\coprod_{t\in T}P^{nA}(\ker(\kappa^{nA}_{\mathcal{E}})_{t}),F_p) &= \prod_{t\in T} \uHom_{\ttBan^{nA}_{R}}(P^{nA}(\ker(\kappa^{nA}_{\mathcal{E}})_{t}),F_p)  \\ & = \prod_{t\in T} \uHom_{\ttBan^{A}_{R}}(P^{A}(\ker(\kappa^{nA}_{\mathcal{E}})_{t}),F_p) \\ & = \uHom_{\ttBan^{A}_{R}}(\coprod_{t\in T}P^{A}(\ker(\kappa^{nA}_{\mathcal{E}})_{t}),F_p)
\end{split}
\end{equation}
and
\begin{equation}\label{IndNAd}\begin{split}\uHom_{\ttInd(\ttBan^{A}_{R})}(\mathcal{E}, \iota(F)) & = \ker [\uHom_{\ttInd(\ttBan^{A}_{R})}(P^{A}(\mathcal{E}),F) \to \uHom_{\ttInd(\ttBan^{A}_{R})}(P^{A}(\ker(\kappa^{A}_{\mathcal{E}})),F)] \\ & =  \ker [\lim_{J \subset I,  |J|<\infty} \colim_{k \in K} \uHom_{\ttBan^{A}_{R}}(\coprod_{j\in J}P^{A}(E_j),F_{k}) \\ & \to \lim_{T \subset L, |T|<\infty} \colim_{p\in P} \uHom_{\ttBan^{A}_{R}}(\coprod_{t\in T}P^{A}(\ker(\kappa^{A}_{\mathcal{E}}))_{t},F_p)] \\
\end{split}
\end{equation}
and so we can conclude from equations (\ref{IndNAa}), (\ref{IndNAb}), (\ref{IndNAc}) and (\ref{IndNAd}) that
\begin{equation}\label{eqn:InnerAdjoint}\uHom_{\ttInd(\ttBan^{nA}_{R})}(\pi(\mathcal{E}), F) =\uHom_{\ttInd(\ttBan^{A}_{R})}(\mathcal{E}, \iota(F))
\end{equation}
and so of course $\Hom_{\ttInd(\ttBan^{nA}_{R})}(\pi(\mathcal{E}), F) =\Hom_{\ttInd(\ttBan^{A}_{R})}(\mathcal{E}, \iota(F)).$
Using (\ref{eqn:InnerAdjoint}) we can conclude that $\pi$ interchanges the monoidal structures:
\begin{equation}\label{eqn:piIndInt}
\begin{split}\uHom_{\ttInd(\ttBan^{nA}_{R})}(\pi(U \ootimes_{} V),W) & = \Hom_{\ttInd(\ttBan^{A}_{R})}(U \ootimes_{} V,\iota(W)) \\ & = \Hom_{\ttInd(\ttBan^{A}_{R})}(U ,\uHom_{\ttInd(\ttBan^{A}_{R})}( V,\iota(W))) \\
& =  \Hom_{\ttInd(\ttBan^{A}_{R})}(U ,\iota(\uHom_{\ttInd(\ttBan^{nA}_{R})}(\pi( V),W))) \\
&=  \Hom_{\ttInd(\ttBan^{nA}_{R})}(\pi(U) ,\uHom_{\ttInd(\ttBan^{nA}_{R})}(\pi( V),W)) 
\\ &=  \Hom_{\ttInd(\ttBan^{nA}_{R})}(\pi(U) \ootimes_{} \pi( V),W).
\end{split}
\end{equation}
By the Yoneda lemma for the opposite category of $\ttInd(\ttBan^{nA}_{R})$ we conclude using (\ref{eqn:piIndInt}) for all $W$ that the natural morphism 
\[\pi(U) \ootimes_{\ttInd(\ttBan^{nA}_{R})} \pi( V) \to \pi(U \ootimes_{\ttInd(\ttBan^{A}_{R})} V)
\]
is an isomorphism.

\ \hfill $\Box$

\bibliographystyle{amsalpha}

\end{document}